\documentclass[12pt]{article}

\usepackage[left=1.5in, right=1in, top=1in, bottom=1in, dvips,pdftex]{geometry}
\usepackage{color, setspace, graphicx, verbatim, mathrsfs, amsmath, amsthm, amssymb, amsfonts,  multirow, moreverb}
\newtheorem{defn}{Definition}[section]

\newtheorem{lemma}[defn]{Lemma}
\newtheorem{thm}[defn]{Theorem}
\newtheorem{cor}[defn]{Corollary}

\theoremstyle{remark}

\numberwithin{equation}{section}
 \numberwithin{figure}{subsection}

 \DeclareMathOperator{\airy}{Ai}
\newcommand{\rchi}{\raisebox{.4ex}{$\chi$}}
\def\ra{\rightarrow}
\def\iy{\infty}

\def\be{\begin{equation}}
\def\ee{\end{equation}}
\def\ov{\over}

\newcommand{\bE}{\mathbb{E}}
\newcommand{\bP}{\mathbb{P}}
\newcommand{\bR}{\mathbb{R}}

\begin{document}

\title{\textbf{Edgeworth Expansion of the Largest Eigenvalue Distribution Function
 of GUE and LUE }}

\author{Leonard N.~Choup \\ Department of Mathematics \\ University of California\\
 Davis, CA 95616, USA \\
 email:  \texttt{choup@math.ucdavis.edu}
}
\maketitle

\begin{center} \textbf{Abstract} \end{center}
\begin{small}
We derive expansions of the Hermite and
Laguerre kernels at the edge
of the spectrum of the finite $n$ Gaussian Unitary Ensemble (GUE$_n$) and the finite $n$
Laguerre Unitary Ensemble (LUE$_n$), respectively.   Using these large $n$
kernel expansions, we
prove an Edgeworth type theorem for the largest eigenvalue distribution function
of GUE$_n$ and LUE$_n$.  In our  Edgeworth  expansion,
the  correction terms  are expressed
in terms of the same Painlev\'e II function  appearing in the leading term, i.e. in
the Tracy-Widom distribution.  We conclude with a brief discussion of the universality
of these results.
\end{small}

\section{Introduction}\label{intro}
The limiting distribution function for the  largest eigenvalues in
orthogonal, unitary and symplectic random matrix ensembles have
found many applications outside their initial discovery in random
matrix theory, see, for example,  \cite{Deif3, Dien1, Trac8} for
recent reviews. In  these applications it is important to have
correction terms to the  limiting distribution.  For example, in
statistics \cite{John1} the sample size is always finite; and to
assess quantitatively  the range of  validity of limit laws,  one
needs  finite $n$ correction terms. In classical probability, a
similar issue arises in the
 application of the Central Limit Theorem (CLT) to finite $n$ problems.  Here the
 two main results are the
Berry-Esseen theorem  and
the Edgeworth expansion \cite{Fell2}.

Recall   if $S_n$ is a sum of  i.i.d.\  random variables $X_j$,
each with  mean $\mu$ and variance $\sigma^2$,
 that the distribution $F_n$ of the normalized
random variable $(S_n - n\mu)/(\sigma \sqrt{n})$ satisfies the
Edgeworth expansion\footnote{We
assume, of course,  the moments $\bE(X_j^k)$, $k=3,\ldots, r$,  exist; and  as well,
the condition $\lim_{\vert\zeta\vert\ra\iy}\sup\vert\varphi(\zeta)\vert <\iy$
where $\varphi$ is the characteristic function of $X_j$, see \cite{Fell2}.}
\begin{equation}\label{edgeworth}
 F_{n}(x)-\Phi(x) = \phi(x)\sum_{j=3}^{r} n^{-\frac{1}{2} j
+1}R_{j}(x)+o(n^{-\frac{1}{2}r+1})
\end{equation}
uniformly in $x$.
Here $\Phi$ is the standard normal distribution
with density $\phi$, and  $R_{j}$ are  polynomials depending only on
$\bE(X_j^k)$ but not on $n$ and $r$ (or the underlying distribution
of the $X_j$).

Introduce
\begin{equation}
F^{G,L} _{n,2}(t)= \bP_{G,L}(\lambda_{\textrm{max}}^{G,L} \leq t)
\end{equation}
where $\lambda_{\textrm{max}}^{G,L}$ is the largest eigenvalue in
 GUE$_n$ or  LUE$_n$, respectively.  (When the context is clear, we
 often drop the  $G$ or $L$.)
To obtain a nontrivial limit theorem, we must, as is well known,  define
 normalized random variables $\hat{\lambda}_{\textrm{max}}^{G,L}$.  We find it useful to ``fine tune'' our normalization (see also \cite{John1}),
 \begin{eqnarray}
 \hat{\lambda}_{\textrm{max}}^G&:=& { \lambda_{\textrm{max}}^G-\left(2(n+c_G)\right)^{1/2}\ov
 2^{-1/2} n^{-1/6}}\, ,\label{GUENorm} \\
 \hat{\lambda}_{\textrm{max}}^L&:=& {\lambda_{\textrm{max}}^L-4(n+c_L)-2\alpha\ov
 2 (2n)^{1/3}}\label{LUENorm}
 \end{eqnarray}
 where $c_{G,L}$ are constants (to be chosen later), and $\alpha$ is the parameter
 appearing in LUE$_n$.  (That is, the parameter $\alpha$ appearing
 in the Laguerre polynomials $L_n^{\alpha}$.)  Then   $\hat{\lambda}_{\textrm{max}}^{G,L}$
 converge in distribution to GUE Tracy-Widom (commonly
 denoted  $F_2$).
 In this paper we initiate the study of
  Edgeworth expansions for both GUE$_n$ and
  LUE$_n$; that is,  we find the analogue of (\ref{edgeworth})
 for $F_n^{G,L}-F_2$.  We now state our main results.

Our first result is an extension of the  Plancherel-Rotach theorem for the Laguerre poynomials
$L_n^{\alpha}$.  We set
\begin{equation}
 \xi =(4n +2\alpha +2c)^{\frac{1}{2}} +
\frac{X}{2^{\frac{2}{3}}n^{\frac{1}{6}}} \quad  \textmd{where} \;
X \; \textmd{and} \; c \;\textmd{are} \; \textmd{bounded},
\end{equation}
and denote by $\airy$  the Airy function  (see, e.g.,  \cite{Olve1}).

\begin{thm}\label{plancherel for Laguerre}
 For $ \alpha > -1$ we have as $n \rightarrow \infty$
\begin{displaymath}
e^{-\xi^2/2}
L_{n}^{\alpha}(\xi^{2})=(-1)^{n}2^{-\alpha-\frac{1}{3}}n^{-\frac{1}{3}}
\quad \left\{  \airy(X) +
 \frac{(c-1)}{2^{\frac{1}{3}}} \airy^{\prime}(X) n^{-\frac{1}{3}} +
 \right.
 \end{displaymath}
 \begin{displaymath}
\left[ \frac{2-10c+5c^2 -5 \alpha}{10 \cdot 2^{\frac{2}{3}} }\, X
\airy(X) + \frac{X^2}{20 \cdot 2^{\frac{2}{3}}}
 \airy^{\prime}(X)\right]n^{-\frac{2}{3}} +
 \end{displaymath}
 \begin{displaymath}
\left[ (\frac{5 \alpha -15 c \alpha + 2 c^3 -15 c^2 -56 c -6}{60}
+ \frac{c-1}{40} X^3 ) \airy(X) + \right.
\end{displaymath}
 \begin{displaymath}
 \left. \left. \frac{(c-1)(5(c-2)c - 3(2+ 5
\alpha))
 }{60} \, X \airy^{\prime}(X) \right]
n^{-1}+ O(n^{-\frac{4 }{3}}) \airy(X) \right
 \}
\end{displaymath}
\end{thm}

From Theorem \ref{plancherel for Laguerre} we derive expansions for both
the Hermite and Laguerre kernels.  Recall that if\hspace{.2ex}\footnote{Here $H_n(x)$ are the
Hermite polynomials of degree $n$.}
\begin{equation*}
 \varphi_{n}(x)= {1\over (2^{n}n! \sqrt{\pi})^{1/2}} \, H_n(x)\,  e^{-x^2/2} \quad \textrm{and}
\quad
\phi_{n}^{\alpha}(x)=x^{\alpha/2}\, e^{-x/2}\, L_{n}^{\alpha}(x),
\end{equation*}
then
the Hermite kernel is
\begin{equation}\label{kernel-her}
K_{n}(x,y)=\sum_{k=0}^{n-1}\varphi_{k}(x)\varphi_{k}(y)
\end{equation}
and the Laguerre kernel is
\begin{equation}\label{kernel-lag}
K_{n}^{\alpha}(x,y)=\sum_{k=0}^{n-1}\frac{\phi_{k}^{\alpha}(x)\phi_{k}^{\alpha}(y)}{\Gamma(k+1)
\Gamma(\alpha +k+1)}\, .
\end{equation}
Finally, the Airy kernel is
\begin{equation}
K_{\airy }(x,y)=
\frac{\airy(x) \airy^{'}(y) - \airy(y) \airy^{'}(x)}{x-y}=\int_0^{\iy} \airy(x+z)\airy(y+z)\, dz.\label{airyKernel}
\end{equation}
Using Theorem \ref{plancherel for Laguerre} we prove
\begin{thm}\label{Hermite kernel}
For $x  =  \left(2(n+c_{G})\right)^{\frac{1}{2}}+2^{-\frac{1}{2}}n^{-\frac{1}{6}}X$ and
$ y =\left(2(n+c_{G})\right)^{\frac{1}{2}}+2^{-\frac{1}{2}}n^{-\frac{1}{6}}Y$
with $X,\, Y $ and $c_{G}$ bounded,
\begin{eqnarray*}
K_{n}(x,y)\,dx =  \left\{  K_{\airy}(X,Y) -c_{G}\airy(X)\airy(Y)
n^{-\frac{1}{3}} +   \right.
\end{eqnarray*}
\begin{eqnarray*}
\frac{1}{20}\left[(X + Y)\airy^{\prime}(X)\airy ^{\prime}(Y) -
(X^2+XY+Y^2)\airy(X)\airy(Y) + \right.
\end{eqnarray*}
\begin{equation}\label{hermite kernel}
\left. \left. \frac{-20c_{G}^2 +3 }{2}(\airy ^{\prime}(X)\airy(Y) +
\airy(X)\airy ^{\prime}(Y)) \right]n^{-\frac{2}{3}}  +O(n^{-1}) E(X,Y)
\right\} \, dX.
\end{equation}
 The error term, $ E(X,Y)$, is the kernel of an integral operator
 on $L^{2}(J)$ which is trace class for any Borel subset $J$ of the reals
that is bounded away from minus infinity.
\end{thm}
For the Laguerre kernel  we prove
\begin{thm} \label{Laguerre kernel}
For $x =4(n+c_{L}) +2\alpha  +2(2n)^{\frac{1}{3}}X$
and  $y =4(n+c_{L}) +2\alpha  +2(2n)^{\frac{1}{3}}Y \;$  with
 $X$,  $Y$ and $c_{L}$ bounded,
 \begin{eqnarray*}
 K_{n}^{\alpha}(x,y)\, dx= \left\{  K_{\airy}(X,Y)
-2^{\frac{2}{3}}c_{L}\airy(X)\airy(Y) n^{-\frac{1}{3}} + \right.
\end{eqnarray*}
\begin{eqnarray*}
\frac{2^{\frac{1}{3}}}{ 10} \left[ \quad (X^2 +X\,Y
+Y^2)\airy(X)\airy(Y) -(X+Y)\airy^{\prime}(X)\airy^{\prime}(Y) - \right.
\end{eqnarray*}
\begin{equation}\label{laguerre kernel}
\left. \left. (10c_{L}^2-1)(\airy(X)\airy^{\prime}(Y)
+\airy^{\prime}(X)\airy(Y)) \right]n^{-\frac{2}{3}} + O(n^{-1})F(X,Y)
\right\} \,dX.
\end{equation}
The error term $F(X,Y)$ is the kernel of an integral operator on
$L^{2}(J)$ which is trace class for any Borel subset $J$ of the reals
which is bounded away from minus infinity.
\end{thm}
To state our main theorem, we need a number of definitions.  First define
the constants
\begin{equation*}
 a_{c_{G},2}^{G}=c_{G},\quad
 a_{c_{L},2}^{L}=2^{\frac{2}{3}}c_{L}, \quad
 b_{2}^{G}=-\frac{1}{20}, \quad
 b_{2}^{L}=\frac{2^{\frac{1}{3}}}{10}.
\end{equation*}
Following the notations of Tracy and Widom
\cite{Trac7}, we set
\begin{equation*}
E_{c_{G},2}^{G}(s)=2w_{1}-3u_{2}
+(-20c_{G}^{2}+3)v_{0}+u_{1}v_{0}-u_{0}v_{1}+u_{0}v_{0}^{2}-u_{0}^{2}w_{0},
\end{equation*}
\begin{eqnarray*}
E_{c_{L},2}^{L}(s)=2w_{1}-3u_{2}
+(20c_{L}^{2}-2)v_{0}+u_{1}v_{0}-u_{0}v_{1}+u_{0}v_{0}^{2}-u_{0}^{2}w_{0},
\end{eqnarray*}
where
\begin{equation*}
 u_{i}:=u_{i}(s)= \int_{s}^{\iy} q(x)x^{i} \airy(x)\,dx, \quad
v_{i}:=v_{i}(s)= \int_{s}^{\iy} q(x)x^{i} \airy^{\prime}(x)\, dx \quad
\textrm{and}
\end{equation*}
\begin{equation*}
 w_{i}:=w_{i}(s)=\int_{s}^{\iy} q'(x)x^{i}
\airy^{\prime}(x)\, dx \;+\; u_{0}(s)v_{i}(s)
\end{equation*}
with $q$ the solution of the Painlev\'{e} II equation, $q'' =s q +2
q^{3}$,  subject to the boundary condition $q(s) \sim \airy(s)$ as
$s \rightarrow \iy$.
Finally, the GUE Tracy-Widom distribution is
\begin{equation}
F_2(s)=\det\left(I-K_{\airy}\rchi_{(s,\iy)}\right)=\exp\left(-\int_s^{\iy} (x-s) q(x)^2\, dx\right)
\label{TW2}
\end{equation}
where $K_{\airy}$ is the operator with Airy kernel (\ref{airyKernel}) and $\rchi_{(s,\iy)}$ is
the indicator function of the interval $(s,\iy)$.
With now can state our main result.

\begin{thm}\label{Edgeworth expansion}
We set
\begin{equation}\label{scale for GUE}
t  =  (2(n+c_{G})
)^{\frac{1}{2}}+2^{-\frac{1}{2}}n^{-\frac{1}{6}}\,s\>\>\> \textrm{for GUE}_n
\end{equation}
and
\begin{equation}\label{scale for LUE}
 t =4(n+c_{L}) +2\alpha  +2(2n)^{\frac{1}{3}}\,s\>\>\>\textrm{for LUE}_n\, .
\end{equation}
Then as $n\ra\iy$
\begin{equation}
F_{n,2}^{G,L}(t)= F_{2}(s)\{ 1 + a_{c_{G,L},2}^{G,L} \,
u_{0}(s) \, n^{-\frac{1}{3}}  +
b_{2}^{G,L}E_{c_{G,L},2}^{G,L}(s)\, n^{-\frac{2}{3}}\} + O(n^{-1})
\end{equation}
uniformly in $s$.
If in addition,
\begin{equation*}
c_{G}^{2}+c_{L}^{2}=\frac{1}{4},\>\> \textrm{then}\>\>
E_{{c_{G}},2}^{G}(s)=E_{{c_{L}},2}^{L}(s) =E_{c,2}(s),\>\>
\textrm{and}
\end{equation*}
\begin{equation}\label{universal}
F_{n,2}^{G,L}(t)= F_{2}(s)\{ 1 + a_{c_{G,L},2}^{G,L} \,
u_{0}(s) \, n^{-\frac{1}{3}}  + b_{2}^{G,L}E_{c,2}(s)\,
n^{-\frac{2}{3}}\} + O(n^{-1}).
\end{equation}
\end{thm}
Note that the $n^{-\frac{1}{2}}$ correction term in the Edgeworth
expansion for the CLT is universal in the sense that only the
constant factor depends on the underlying  distribution.   We see in (\ref{universal})
a similar universality, and we  conjecture that
this universality extends to a wider class of unitary ensembles.

 In \S 2 we  derive Theorems
\ref{Hermite kernel} and \ref{Laguerre kernel} from
Theorem \ref{plancherel for Laguerre}. In \S3 we follow
\cite{Trac7} to prove Theorem \ref{Edgeworth expansion}. The proof of
Theorem \ref{plancherel for Laguerre} will be given in the
Appendices together with some facts needed to prove the last
Theorem.

\section{Correction terms for the Hermite and Laguerre kernel at
the edge of the spectrum } To simplify notations we will use
matrix ensembles of $(n+1)\times( n+1)$ matrices throughout this
section and part of the next section. After the fine tuning of the
variables in \S\ref{fineTuning}, we will use ensembles of $n \times n$
matrices.

\subsection{Hermite case}
We have the following representation of the Hermite kernel from
the Christoffel-Darboux formula.
\begin{equation}\label{HermiteK}
 K_{n+1}(x,y)  =
\sqrt{\frac{n+1}{2}} \quad \frac{\varphi_{n+1}(x)\varphi_{n}(y)
 -  \varphi_{n+1}(y)\varphi_{n}(x)}{x-y}
\end{equation}
As mentioned in \S\ref{intro}, Theorem \ref{Hermite kernel} is
a corollary of  Theorem \ref{plancherel for
Laguerre}. We recall the relation between the Hermite and Laguerre polynomials\\
\[ H_{2n}(x)=(-1)^{n}2^{2n}n!L_{n}^{-\frac{1}{2}}(x^{2}), \>\>
H_{2n+1}(x)=(-1)^{n}2^{2n+1}n!xL_{n}^{\frac{1}{2}}(x^{2})\]
so as to estimate the right side of \eqref{HermiteK}. We
will also assume\footnote{For $n$ odd the same analysis can be
carried out to produce the same result.} without lost of
generalities that $n=2k$ is even.
Using the symmetry of the Hermite kernel in \eqref{HermiteK} we
only need to find an expansion of
\begin{equation}\label{eq 2}
\varphi_{n+1}(x)\varphi_{n}(y) =
 \frac{ 2^{4k+1}(k!)^2
xe^{-\frac{x^{2}}{2}}L_{k}^{\frac{1}{2}}(x^{2})
e^{-\frac{y^{2}}{2}}L_{k}^{-\frac{1}{2}}(y^{2}) }
{(2^{4k+1}(2k+1)! (2k)!\pi)^{\frac{1}{2}}}
\end{equation}
since the other term follows by interchanging $x$ and $y$. The
Laguerre polynomial of argument  $y$ in \eqref{eq 2} has parameter
$\alpha=-\frac{1}{2}$; thus,  in order to apply Theorem \ref{plancherel for
Laguerre}, we write $y$ in the form $y = (4k-1
+2c^{'})^{\frac{1}{2}} + \frac{ Y}{2^{\frac{2}{3}}
k^{\frac{1}{6}}}$ which  corresponds to $c^{'}=c +1$.
Next we use  Stirling's  formula  to estimate
\begin{equation*}
 \quad \frac{ 2^{4k+1}(k!)^2 x } {(2^{4k+1}(2k+1)!
(2k)!\pi)^{\frac{1}{2}}}(-1)^{k} 2^{-\frac{1}{2}-\frac{1}{3}}
k^{-\frac{1}{3}} (-1)^{k}
2^{\frac{1}{2}-\frac{1}{3}}k^{-\frac{1}{3}}\, .
\end{equation*}
The last factors in the left hand side are the constant factor in
Theorem \eqref{plancherel for Laguerre} for the Laguerre functions in $x$
and $y$ respectively. This expression is
\begin{equation*}
 \quad 2^{\frac{1}{3}} k^{-\frac{1}{6}}\left( 1+ \frac{X}{2^{\frac{5}{3}}
k^{\frac{2}{3}}} + \frac{c}{4k} +O(k^{-\frac{5}{3}})\right) .
\end{equation*}
This times,  the constant\footnote{Recall that in this section the
notation $K_{n}$  stands for an ensemble of $(n+1)\times (n+1)$
matrices.} $\sqrt{\frac{n+1}{2}}$ from \eqref{HermiteK},  gives
\begin{equation*}
2^{\frac{1}{3}} k^{\frac{1}{3}}( 1+ \frac{X}{2^{\frac{5}{3}}
k^{\frac{2}{3}}} + \frac{c+1}{4k} +O(k^{-\frac{5}{3}}) ) .
\end{equation*}

We substitute all these into \eqref{eq 2},  and then interchange $x$ and
$y$ to have the second term $\varphi_{n+1}(y)\varphi_{n}(x)$.
Finally with the help of {\em Mathematica} we derive the following
version of Theorem \ref{Hermite kernel}.

\begin{equation}\label{changement de varriables for H}
\textrm{For} \>\>  x =  (2n+1 +2c)^{\frac{1}{2}} +
\frac{ X}{2^{\frac{1}{2}} n^{\frac{1}{6}}}\>\>  \textrm{and}\>\>
 y =  (2n+1 +2c)^{\frac{1}{2}} + \frac{Y}{2^{\frac{1}{2}} n^{\frac{1}{6}}}\, ,
\end{equation}
\begin{eqnarray*}
K_{n+1}(x,y)dx = \left \{  K_{\airy}(X,Y) +\frac{1-2c}{2}
\airy(X)\airy(Y) n^{-\frac{1}{3}} + \left[\frac{X +
Y}{20}\airy^{'}(X)\airy^{'}(Y) -\right. \right.
\end{eqnarray*}
\begin{eqnarray*}
\left. \frac{X^2+XY+Y^2}{20}\airy(X)\airy(Y) + \frac{-10c^2 +10c
-1}{20}(\airy^{'}(X)\airy(Y) + \airy(X)\airy^{'}(Y)) \right]n^{-\frac{2}{3}}
\end{eqnarray*}
\begin{equation}\label{hermitekernel}
\left. + \,O(\,\frac{1}{n}\,) \;E(X,Y)  \right\} dX
\end{equation}
Using the
symmetry between  the $X$ and  $Y$  terms from
\eqref{HermiteK}, this error term can be expressed as a finite sum
\begin{eqnarray*}
E(X,Y)&=& \sum_{j}P_{j}(X,Y) \airy(X)\airy(Y) +
\sum_{j_{1}}Q_{j_{1}}(X,Y) \airy^{'}(X)\airy(Y) +  \\
& &\sum_{j_{2}}Q_{j_{2}}(X,Y) \airy(X)\airy^{'}(Y) +
\sum_{j}R_{j}(X,Y) \airy^{'}(X)\airy^{'}(Y) +
\\
& &\sum_{j,k}a_{j,k}\frac{X^{j}Y^{k}\airy(X)\airy^{'}(Y)-
X^{k}Y^{j}\airy^{'}(X)\airy(Y)}{X-Y}
\end{eqnarray*}
where all coefficients of the polynomials $P$, $Q$ and $R$, and
the $a_{j,k}$ have a factor of $n^{-\frac{1}{3}k}$, $k
\in \{0,1,2,3 \}$.

The first four terms are kernels of
finite rank operators on any Borel subset $J$ of $\bR$
not including minus infinity;\footnote{This last restriction is due to the
behavior of the Airy function near minus infinity.} and thus, are trace
class. The last term also defines a  trace class operator. This is best seen from
the following result. If we assumed without lost of generalities
that $j\leq k$,  and set  $k-j=s$, then
\begin{eqnarray*}
&
&\frac{X^{j}Y^{k}\airy(X)\airy^{'}(Y)-X^{k}Y^{j}\airy^{'}(X)\airy(Y)}{X-Y}\\
&=&(XY)^{j}
\frac{(Y^{s}-X^{s}+X^{s})\airy(X)\airy^{'}(Y)-X^{s}\airy^{'}(X)\airy(Y)}{X-Y}\\
&=& (XY)^{j}\left( X^{s}K_{\airy}(X,Y)\;+\;
\sum_{i=0}^{s-1}Y^{s-i}X^{i}\airy(X)\airy^{'}(Y) \right).
\end{eqnarray*}
This shows that the error $E(X,Y)$ in \eqref{hermitekernel} is the
kernel of a trace class operator.\\
\subsection{Laguerre case}
Again by the  Christoffel-Darboux formula,
\begin{equation}\label{eq 3}
 K_{n+1}^{\alpha}(x,y)  =
\frac{(n+1)(xy)^{\frac{\alpha}{2}}}{\Gamma(\alpha
+1)\left(\begin{array}{c}
  n+\alpha \\
  n \\\end{array} \right)}
\frac{\phi_{n+1}^{\alpha}(x)\phi_{n}^{\alpha}(y)  -
\phi_{n+1}^{\alpha}(y)\phi_{n}^{\alpha}(x)}{x-y}\, .
\end{equation}

\subsubsection{Asymptotic of $ e^{-\frac{x}{2}}L_{n}^{\alpha}(x)$
at $x=4n +2\alpha +2c + 2(2n)^{\frac{1}{3}}X $}

In order to apply Theorem \ref{plancherel for Laguerre} to the
Laguerre kernel at the edge of the spectrum (corresponding to
$x=4n +2\alpha +2c + 2(2n)^{\frac{1}{3}}X $ for bounded
$c$ and $X$), we need to make a variable change $x=\xi^2$ where
$\xi =\sqrt{4n +2\alpha +2c} +
\;2^{-\frac{2}{3}}n^{-\frac{1}{6}}\; t $. We use these two
expressions to solve for $t$ in terms of $X$ and then substitute
this value for $t$ into Theorem \ref{plancherel for Laguerre} to obtain
the desired asymptotics. But\footnote{See the last footnote in the
proof  of Theorem \ref{plancherel for Laguerre} in the Appendix for
this technical point.} in order to have  accurate asymptotics for
the Laguerre functions at the edge of the spectrum, we will use
the expression of $\xi$ involving $l_{n}=(4n +2\alpha +2c)^{1/2}$.
Thus
\begin{equation*}
x=\xi^2 \quad \textrm{is equivalent to } \quad l_{n}^2
+2(2n)^{\frac{1}{3}}X=l_{n}^2 +(2l_{n})^{\frac{2}{3}}t
+(2l_{n})^{-\frac{2}{3}}t^2 \; .
\end{equation*}
\begin{eqnarray*}
\textrm{This quadratic has solutions} \; t_{\pm}=\frac{1}{2}
(-(2l_{n})^{\frac{4}{3}} \pm (2l_{n})^{\frac{4}{3}}\sqrt{1
+4(2l_{n})^{-2} 2(2n)^{\frac{1}{3}}X})\; .
\end{eqnarray*}
Actually, only the solution with the plus sign is to be taken as
it is the only one bounded when $n$ increases. An expansion of
$t_{+}=t$ leads to:
\begin{eqnarray*}
t=X-\frac{(\alpha +c)}{6n}X -
\frac{X^2}{2^{\frac{8}{3}}n^{\frac{2}{3}}} +O(n^{-\frac{4}{3}})\,.
\end{eqnarray*}

Thus if in Theorem \ref{plancherel for Laguerre} we replace $X$ by
this value of $t$, we obtain, again  with the help of {\em Mathematica},  the
desired expansion.
\begin{lemma}\label{lemma for l_{n}}
 For $\; x =4n +2\alpha +2c + 2(2n)^{\frac{1}{3}}X
 $ and $X$ bounded,
\begin{equation*}
 e^{-\frac{x}{2}}
L_{n}^{\alpha}(x)=(-1)^{n}2^{-\alpha-\frac{1}{3}}n^{-\frac{1}{3}}
\left \{ \airy(X) +
 \frac{(c-1)}{2^{\frac{1}{3}}} \airy^{'}(X) n^{-\frac{1}{3}} +
 \right.
 \end{equation*}
 \begin{equation*}
\left[ \frac{2-10c+5c^2 -5 \alpha}{10 \cdot 2^{\frac{2}{3}} }X
\airy(X) - \frac{2 X^2}{10 \cdot 2^{\frac{2}{3}}}
 \airy^{'}(X)\right]n^{-\frac{2}{3}} +
 \end{equation*}
 \begin{equation*}
\frac{1}{60} \left[ -\left(6+56c +15c^2 -2c^3 +5\alpha(3c-1) -6X^3
+6cX^3 \right) \airy(X)  \right.
\end{equation*}
\begin{equation}\label{eq 4}
 \left. \left. + \left(6+\alpha(5-15c)-6c -15c^2 +5c^3\right) X \airy^{'}(X)\right]
n^{-1}\; +\; O(n^{-4/3}) \, \airy(X)
\right \}
\end{equation}
\end{lemma}

\subsubsection{Asymptotic of $ e^{-\frac{x}{2}}L_{n+1}^{\alpha}(x)$
at $x=4n +2\alpha +2c + 2(2n)^{\frac{1}{3}}X $}

Making use of this last formula, we can derive an asymptotic for\\
$e^{-\frac{x}{2}} L_{n+1}^{\alpha}(x) $ when  $ x =4n +2\alpha +2c
+ 2(2n)^{\frac{1}{3}}X$. Note that the degree of the Laguerre
polynomial is no longer $n$ but $n+1$, so in order to
use Lemma\ref{lemma for l_{n}}, we need to write $x$ in terms of $n+1$ or
 $ x=4(n+1) +2\alpha
+2(c-2) + 2(2(n+1))^{\frac{1}{3}}[X- \frac{1}{3n}X +O(n^{-2})]$.\\
The substitution needed here is $c\rightarrow c-2 \,$, and $ X
\rightarrow X - \frac{1}{3n}X$.  {\em Mathematica}
again gives:
\begin{lemma}\label{lemma for l_{n+1}}
 For $x =4n +2\alpha +2c + 2(2n)^{\frac{1}{3}}X
 $ and $X$ bounded,
\begin{eqnarray*}
 e^{-\frac{x}{2}}
L_{n+1}^{\alpha}(x)=(-1)^{n+1}2^{-\alpha-\frac{1}{3}}(n+1)^{-\frac{1}{3}}
\left\{ \airy(X) +
 \frac{(c-3)}{2^{\frac{1}{3}}} \airy^{'}(X) n^{-\frac{1}{3}} +
 \right.
 \end{eqnarray*}
 \begin{eqnarray*}
\left[ \frac{42-30c+5c^2 -5 \alpha}{10 \cdot 2^{\frac{2}{3}} }X
\airy(X) - \frac{2 X^2}{10 \cdot 2^{\frac{2}{3}}}
 \airy^{'}(X)\right]n^{-\frac{2}{3}} +
 \end{eqnarray*}
 \begin{eqnarray*}
\frac{1}{60}  \left[ \left(30+28c -27c^2 +2c^3 -5\alpha(3c-7)
+18X^3 -6cX^3 \right) \airy(X) \right.
\end{eqnarray*}
\begin{equation}\label{eq 5}
\left. \left. + (-102-5 \alpha(3c-7)+114c -45c^2 +5c^3) X
\airy^{'}(X)\right] n^{-1}\; + \;O(n^{-4/3}) \, \airy(X) \right\}
\end{equation}
\end{lemma}

To complete this subsection we need to estimate
\begin{equation*}
\frac{(n+1)(xy)^{\frac{\alpha}{2}}}{\Gamma(\alpha
+1)\left(\begin{array}{c}
  n+\alpha \\
  n \\
\end{array} \right)}.
\end{equation*}
We have
\begin{eqnarray*}
 (xy)^{\frac{\alpha}{2}}=2^{2\alpha}n^{\alpha}
[1+\frac{\alpha(X+Y)}{2^{\frac{5}{3}}n^{\frac{2}{3}}} +
\frac{\alpha(\alpha +c)}{2n} +O(n^{-\frac{4}{3}})],
\end{eqnarray*}
\begin{eqnarray*}
\Gamma(\alpha +1)\left(\begin{array}{c}
  n+\alpha \\
  n \\
\end{array} \right) =n^{\alpha}[1+\frac{\alpha^2 +\alpha}{2n}
+O(n^{-2})]
\end{eqnarray*}
and in addition, the product of the constant factor on the right
of \eqref{eq 4} and \eqref{eq 5} is
\begin{eqnarray*}
(-1)^{n+1}2^{-\alpha -\frac{1}{3}}(n+1)^{-\frac{1}{3}}
\cdot(-1)^{n}2^{-\alpha -\frac{1}{3}}n^{-\frac{1}{3}}
=-2^{-2\alpha}n^{-\frac{2}{3}}[1 -\frac{1}{3n} +O(n^{-2})]
\end{eqnarray*}
thus
\begin{eqnarray*}
\frac{(n+1)(xy)^{\frac{\alpha}{2}}}{\Gamma(\alpha
+1)\left(\begin{array}{c}
  n+\alpha \\
  n \\
\end{array} \right)}\cdot
(-1)^{n+1}2^{-\alpha -\frac{1}{3}}(n+1)^{-\frac{1}{3}}
\cdot(-1)^{n}2^{-\alpha -\frac{1}{3}}n^{-\frac{1}{3}}=
\end{eqnarray*}
\begin{eqnarray*}
-2^{-\frac{2}{3}}n^{\frac{1}{3}}[1+\frac{\alpha(X+Y)}
{2^{\frac{5}{3}}n^{\frac{2}{3}}} +\frac{3\alpha(c-1) +4}{6n}
+O(n^{-\frac{4}{3}})]
\end{eqnarray*}
Substituting all these quantities in \eqref{eq 3}, give the
following version of Theorem  \ref{Laguerre kernel}
\begin{equation}\label{changement de varriables for L}
\textrm{For} \>\>  x  = 4n+2\alpha+2c + 2(2
n)^{\frac{1}{3}}X \>\>  \textrm{and} \>\> y =  \
4n+2\alpha+2c + 2(2 n)^{\frac{1}{3}}Y
\end{equation}
\begin{eqnarray*}
K_{n+1}^{\alpha}(x,y)dx = \left \{  K_{\airy}(X,Y)
+\frac{2-c}{2^{\frac{1}{3}}} \, \airy(X)\, \airy(Y)\,
n^{-\frac{1}{3}} +\right.
\end{eqnarray*}
\begin{eqnarray*}
\frac{1}{2^{\frac{2}{3}}\,10}\left[2(X^{2}+X\,Y +Y^{2})\,
\airy(X)\,\airy(Y) - \; 2(X+Y)\, \airy^{'}(X) \, \airy^{'}(Y)
-\right.
\end{eqnarray*}
\begin{equation}\label{laguerrekernel}
\left. \left.  (18-20c+5c^{2})(\airy(X)\, \airy^{'}(Y) +
\airy^{'}(X) \, \airy(Y)) \right]n^{-\frac{2}{3}}
 + \,O(\,\frac{1}{n}\,) \;F(X,Y)  \right\} dX.
\end{equation}
As in the Hermite case, the error term $F(X,Y)$  is a finite sum
of kernels of trace class operators, therefore is a kernel of a
trace class operator on $L^{2}(J)$ for any subset  $J$ of the
reals which is bounded away from minus infinity.

\subsection{Conclusion}
For our order of expansion, we see that both kernels are finite
rank perturbation of the Airy kernel.
In the Laguerre case, the final result does not involve an
explicit presence of the order $\alpha$.

\section{Expansion of the Fredholm determinants  at the edge of the spectrum}

This part will only make use  \eqref{changement de varriables for
H} and \eqref{hermitekernel} to derive the desired result in the
Hermite case and \eqref{changement de varriables for L} and
\eqref{laguerrekernel} in the Laguerre case. Most of the
derivations will follow from the work of Tracy and Widom.

\subsection{Edgeworth expansion of  $F_{n,2}^{G}$ }
Recall that
\begin{equation}\label{eq 6}
F_{n+1,2}^{G}(t)=\textrm{det}(I-K_{n+1})
\end{equation}
where this determinant is the Fredholm determinant of the integral
operator with kernel $K_{n+1}(x,y)$  on $L^{2}(s,\infty)$,
$t$ and $s$ are related by
\begin{equation}\label{hermitescale}
t=(2\, n  +  1+ 2\,c)^{\frac{1}{2}} +
2^{-\frac{1}{2}}\,n^{-\frac{1}{6}}\,s\, .
\end{equation}
In this section we will estimate this determinant. Most of our
derivations involve trace class operators where the Fredholm
determinant is analytic. If in \eqref{eq 6} we use the expression
of $\;K_{n+1}\;$ given by \eqref{hermitekernel}, the continuity of
the determinant (in trace class norm) allows us to pull the error term involving the
kernel $E(X,Y)$ out of the determinant as an $O(n^{-1})$ term. We
therefore have
\begin{eqnarray*}
 \textrm{det}(I-K_{n+1})=\textrm{det}\left(I-\left\{K_{\airy}(X,Y) +\frac{1-2c}{2} \airy(X)\airy(Y)
n^{-\frac{1}{3}} + \right. \right.
\end{eqnarray*}
\begin{eqnarray*}
\left[\frac{X + Y}{20}\airy^{'}(X)\airy^{'}(Y)
-\frac{X^2+XY+Y^2}{20}\airy(X)\airy(Y) + \right.
\end{eqnarray*}
\begin{equation}\label{eq4}
\left. \left. \left.\frac{-10c^2 +10c -1}{20}(\airy^{'}(X)\airy(Y)
+ \airy(X)\airy^{'}(Y))
\right]n^{-\frac{2}{3}}\right\}\right)+O(n^{-1}).
\end{equation}
Note we are using the obvious notation of writing the kernel for the operator appearing
in the determinant.  We continue to employ this notation below.
If we factor out $\textrm{det}(I-K_{\airy}(X,Y))$ and set
\[ K_{n+1}(x,y)\,dx -K_{\airy}(X,Y)\, dX:=L(X,Y)\,dX +O(n^{-1}) E(X,Y)\,dX, \]
the multiplicative property of the determinant gives
\begin{equation}\label{factor}
 \textrm{det}(I-K_{n+1}(x,y))=\textrm{det}(I-K_{\airy}(X,Y))
 \textrm{det}(I-(I-K_{\airy}(X,Y))^{-1}L(X,Y)) +O(n^{-1}).
\end{equation}
The first factor on the right is known, so we need to estimate the
second factor. We will express this factor as a finite sum of rank
one operators
\begin{equation}\label{eq 7}
(I-K_{\airy}(X,Y))^{-1}L(X,Y)=\sum_{i=1}^{k}\alpha_{i}(X)\beta_{i}(Y)
\end{equation}
and use the well known formula
\begin{equation}\label{eq 8}
\det\left(I-\sum_{i=1}^{k}\alpha_{i}(X)\beta_{i}(Y)\right)=
\det\left(\delta_{ij} -(\alpha_{i},\beta_{j})\right)_{1\leq i,j\leq k}
\end{equation}
to transform the problem to one involving estimations of
inner products. Let $(\cdot,\cdot)$ denote the inner product in
$L^{2}(s,\infty)$. A good reference for much of the results that follows is \cite{Trac7}.
We have the following representation.
\begin{eqnarray*}
 (I-K_{\airy}(X,Y))^{-1}L(X,Y)=\frac{1-2c}{2} (I-K_{\airy}(X,Y))^{-1}\airy(X)\airy(Y)
n^{-\frac{1}{3}} +
\end{eqnarray*}
\begin{eqnarray*}
\frac{1}{20}\left[(I-K_{\airy}(X,Y))^{-1}X
\airy^{'}(X)\airy^{'}(Y) + (I-K_{\airy}(X,Y))^{-1}Y
\airy^{'}(X)\airy^{'}(Y)- \right.
\end{eqnarray*}
\begin{eqnarray*}
(I-K_{\airy}(X,Y))^{-1}X^2\airy(X)\airy(Y) -
(I-K_{\airy}(X,Y))^{-1}X Y \airy(X)\airy(Y) -
\end{eqnarray*}
\begin{eqnarray*}
(I-K_{\airy}(X,Y))^{-1}Y^2\airy(X)\airy(Y)+ (-10c^2 +10c
-1)(I-K_{\airy}(X,Y))^{-1}\airy^{'}(X)\airy(Y)
\end{eqnarray*}
\begin{eqnarray*}
 \left.+ (-10c^2
+10c -1)(I-K_{\airy}(X,Y))^{-1}\airy(X)\airy^{'}(Y)
\right]n^{-\frac{2}{3}} \,+\,O(n^{-1})
\end{eqnarray*}

In the notations above, we think of all the quantities involved as
kernels of integral operators, the analogues of those in
\cite{Trac7}. We therefore have:
\begin{eqnarray*}
\frac{1-2c}{2}(I-K_{\airy}(X,Y))^{-1}\airy(X)\airy(Y)n^{-\frac{1}{3}}=
\frac{1-2c}{2}Q(X)\airy(Y)n^{-\frac{1}{3}}:=
\alpha_{1}(X)\beta_{1}(Y)
\end{eqnarray*}
where $Q(X) $ is the action of the integral operator with kernel
$(I-K_{\airy})^{-1}$ on $\airy(X)$. In the same way we have,
\begin{eqnarray*}
\frac{1}{20}(I-K_{\airy}(X,Y))^{-1}X
\airy^{'}(X)\airy^{'}(Y)n^{-\frac{2}{3}}=
\frac{1}{20}((I-K_{\airy})^{-1}f)(X)
\airy^{'}(Y)n^{-\frac{2}{3}}\\:= \alpha_{2}(X)\beta_{2}(Y)
\end{eqnarray*}
where $f(Z)=Z\airy'(Z)$ and,
\begin{eqnarray*}
\frac{1}{20}(I-K_{\airy}(X,Y))^{-1}\airy^{'}(X)Y
\airy^{'}(Y)n^{-\frac{2}{3}}=
\frac{1}{20}P(X)f(Y)n^{-\frac{2}{3}}\\:= \alpha_{3}(X)\beta_{3}(Y)
\end{eqnarray*}
where $P(X) $ is the action of the integral operator with kernel
$(I-K_{\airy})^{-1}$ acting on $\airy^{'}(X)$.
\begin{eqnarray*}
-\frac{1}{20}(I-K_{\airy}(X,Y))^{-1}X^{2}
\airy(X)\airy(Y)n^{-\frac{2}{3}}=
-\frac{1}{20}((I-K_{Airy})^{-1}g)(x) \airy(y)n^{-\frac{2}{3}}\\:=
\alpha_{4}(X)\beta_{4}(Y)
\end{eqnarray*}
where $g(Z)=Z^{2}\airy(Z)$ ,
\begin{eqnarray*}
-\frac{1}{20}(I-K_{\airy}(X,Y))^{-1}X Y
\airy(X)\airy(Y)n^{-\frac{2}{3}}=
-\frac{1}{20}((I-K_{\airy})^{-1}h)(X) h(Y)n^{-\frac{2}{3}}\\:=
\alpha_{5}(X)\beta_{5}(Y)
\end{eqnarray*}
where $h(Z)=Z \airy(Z)$ ,
\begin{eqnarray*}
-\frac{1}{20}(I-K_{\airy}(X,Y))^{-1}\airy(X)Y^{2}
\airy(Y)n^{-\frac{2}{3}}= -\frac{1}{20}Q(X)g(Y)n^{-\frac{2}{3}}:=
\alpha_{6}(X)\beta_{6}(Y)
\end{eqnarray*}
\begin{eqnarray*}
C\,(I-K_{\airy}(X,Y))^{-1}\airy^{'}(X) \airy(Y)n^{-\frac{2}{3}}= C
P(X)\airy(Y)n^{-\frac{2}{3}}:= \alpha_{7}(X)\beta_{7}(Y) \, ,\quad
\textrm{and}
\end{eqnarray*}
\begin{eqnarray*}
C\,(I-K_{\airy}(X,Y))^{-1}\airy(X) \airy^{'}(Y)n^{-\frac{2}{3}}= C
Q(X)\airy^{'}(Y)n^{-\frac{2}{3}}:= \alpha_{8}(X)\beta_{8}(Y)
\end{eqnarray*}
where $C=(-10c^2+10c-1)/20$.

If we set $(\alpha_{i},\beta_{j})=a_{ij}n^{-\frac{2}{3}}$ for
$j\neq 1$ and
$(\alpha_{i},\beta_{1})=a_{i1}n^{-\frac{1}{3}}$ , expanding
\eqref{eq 8} with respect to the first row leads to the following
expression.
\begin{eqnarray*}
\textrm{det}(\delta_{ij} -(\alpha_{i},\beta_{j}))_{ i,j=1
}^{8} &=& (1- a_{11}n^{-\frac{1}{3}})\left|%
\begin{array}{cccc}
 1-  a_{22}n^{-\frac{2}{3}}& -a_{23}n^{-\frac{2}{3}} & \cdots & -a_{28}n^{-\frac{2}{3}} \\
  -a_{32}n^{-\frac{2}{3}} & 1-a_{33}n^{-\frac{2}{3}} & \cdots & -a_{38}n^{-\frac{2}{3}} \\
  \cdot & \cdot & \cdot & \cdot \\
  -a_{82}n^{-\frac{2}{3}} & -a_{83}n^{-\frac{2}{3}} & \cdots & 1-a_{88}n^{-\frac{2}{3}} \\
\end{array}%
\right| \\
&+&
\sum_{k=2}^{8}(-1)^{k}a_{1k}n^{-\frac{2}{3}}\textrm{det}(C_{1},
\cdots,\hat{C_{k}},\cdots, C_{8}) \\
&=& (1- a_{11}n^{-\frac{1}{3}})\left|%
\begin{array}{cccc}
 1-  a_{22}n^{-\frac{2}{3}}& -a_{23}n^{-\frac{2}{3}} & \cdots & -a_{28}n^{-\frac{2}{3}} \\
  -a_{32}n^{-\frac{2}{3}} & 1-a_{33}n^{-\frac{2}{3}} & \cdots & -a_{38}n^{-\frac{2}{3}} \\
  \cdot & \cdot & \cdot & \cdot \\
  -a_{82}n^{-\frac{2}{3}} & -a_{83}n^{-\frac{2}{3}} & \cdots & 1-a_{88}n^{-\frac{2}{3}} \\
\end{array}%
\right| \\
&+& O(n^{-1}).
\end{eqnarray*}
We factor out $n^{-\frac{1}{3}}$ from column $C_{1}$ in the last
step. The determinant in the last line is of the same form as the
the original determinant, therefore a similar transformation to
this last determinant leads to the following result.
\begin{eqnarray*}
\det \left( \delta_{ij} \right.&-&\left. ( \alpha_{i},\beta_{j}) \right)_{1\leq i,j \leq 8} =
(1-a_{11}n^{-\frac{1}{3}})\prod_{k=2}^{8}(1-a_{kk}n^{-\frac{2}{3}}) + O(n^{-1})\\
&=&(1-
a_{11}n^{-\frac{1}{3}})\sum_{k=0}^{7}(-1)^{k}n^{-\frac{2k}{3}}
\sum_{i_{1},\cdots,i_{k} ;i_{r}\neq i_{s} \in \{2,\dots,8\}} \quad
\prod_{j=1}^{k}a_{i_{j}i_{j}} + O(n^{-1}) \\
&=& 1-a_{11}n^{-\frac{1}{3}} -
n^{-\frac{2}{3}}\sum_{k=2}^{8}a_{kk}+ O(n^{-1}) = 1 -
\sum_{k=1}^{8}(\alpha_{k},\beta_{k}) + O(n^{-1}).
\end{eqnarray*}
Thus we only need to compute the inner products
$(\alpha_{k},\beta_{k})\;\textrm{for} \;k=1,\cdots,8$.\\ To
simplify notations we will write for example $u_{0}$ instead of
$u_{0}(s)$. We therefore have:
\begin{eqnarray*}
(\alpha_{1},\beta_{1})&=& \frac{1-2c}{2}(Q,\airy)n^{-\frac{1}{3}}
= \frac{1-2c}{2}\, u_{0} \, n^{-\frac{1}{3}}\\
\end{eqnarray*}
\begin{eqnarray*}
(\alpha_{2},\beta_{2})&=&
\frac{n^{-\frac{2}{3}}}{20}((I-K_{\airy})^{-1}
X\airy^{'},\airy{'}) = \frac{n^{-\frac{2}{3}}}{20}( P,X\airy^{'})
=\frac{n^{-\frac{2}{3}}}{20}w_{1}\\
\end{eqnarray*}
\begin{eqnarray*}
(\alpha_{3},\beta_{3})&=& \frac{n^{-\frac{2}{3}}}{20}(
P,X\airy^{'}) =\frac{n^{-\frac{2}{3}}}{20}w_{1}\\
\end{eqnarray*}
\begin{eqnarray*}
(\alpha_{4},\beta_{4})&=&
-\frac{n^{-\frac{2}{3}}}{20}((I-K_{\airy})^{-1} g,\airy) =
-\frac{n^{-\frac{2}{3}}}{20}( Q,X^{2}\airy)
=-\frac{n^{-\frac{2}{3}}}{20}u_{2}\\
\end{eqnarray*}
\begin{eqnarray*}
(\alpha_{5},\beta_{5})&=&
-\frac{n^{-\frac{2}{3}}}{20}((I-K_{\airy})^{-1} h,h)=
-\frac{n^{-\frac{2}{3}}}{20}(Q_{1},X \airy)\\
\end{eqnarray*}
\begin{eqnarray*}
(\alpha_{6},\beta_{6})&=& -\frac{n^{-\frac{2}{3}}}{20}( Q,X^{2}
\airy) =-\frac{n^{-\frac{2}{3}}}{20}u_{2}\\
\end{eqnarray*}
\begin{eqnarray*}
(\alpha_{7},\beta_{7})&=& Cn^{-\frac{2}{3}}(
P,\airy)=Cn^{-\frac{2}{3}}( Q,\airy^{'})=Cn^{-\frac{2}{3}}v_{0}
\quad = \quad (\alpha_{8},\beta_{8}).
\end{eqnarray*}
To estimate $(\alpha_{5},\beta_{5})$, we use equation 2.12 of
\cite{Trac7} which says
\begin{equation*}
Q_{1}(X)=XQ(X)+u_{0}P(X)-v_{0}Q(X), \quad \textrm{to have}
\end{equation*}
\begin{eqnarray*}
(Q_{1}(X),X \airy(X))&=&(X\,Q(X)+u_{0}P(X)-v_{0}Q(X),X \airy
(X))\\&=&u_{2}-v_{0}u_{1}+u_{0}(P(X),X\airy(X))\\
&=&u_{2}-v_{0}u_{1}+u_{0}\widetilde{v}_{1} \\
&=&
u_{2}-v_{0}u_{1}+u_{0}v_{1}-u_{0}v_{0}^{2}+u_{0}^{2}w_{0}.\quad
(\textrm{
 We used  }\, (2.14) \, \textrm{ of } \cite{Trac7})
\end{eqnarray*}
And
\begin{equation*}
(\alpha_{5},\beta_{5})=
-\frac{n^{-\frac{2}{3}}}{20}(u_{2}-u_{1}v_{0}+u_{0}v_{1}-u_{0}v_{0}^{2}
+u_{0}^{2}w_{0}).
\end{equation*}
Substituting this into the formula for the determinant gives
\begin{eqnarray*}
& &\textrm{det}(\delta_{ij} -(\alpha_{i},\beta_{j}))_{1\leq
i,j\leq
8}= 1 - \sum_{k=1}^{8}(\alpha_{k},\beta_{k}) + O(n^{-1})\\
&=& 1-\frac{1-2c}{2}u_{0}n^{-\frac{1}{3}}
-\frac{n^{-\frac{2}{3}}}{20} \left\{ 2w_{1}-3u_{2}+ (-20c^2 +20c
-2)v_{0} +u_{1}v_{0} \right.\\ & &- \left.u_{0}v_{1}
+u_{0}v_{0}^{2}-u_{0}^{2}w_{0} \right\} \quad + \quad O(n^{-1}).
\end{eqnarray*}
If we set
\begin{equation*}
a_{c,2}^{G}(s)=\frac{2c-1}{2} \;, \quad b^{G}_{2}=-\frac{1}{20}
\quad \textrm{and}
\end{equation*}
\begin{equation}\label{second term}
E_{c,2}^{G}(s)=2w_{1}-3u_{2}+ (-20c^2 +20c -2)v_{0}  +
u_{1}v_{0}-u_{0}v_{1} +u_{0}v_{0}^{2}-u_{0}^{2}w_{0},
\end{equation}
equation \eqref{factor} gives
\begin{lemma}\label{edgeworth for hermite}
For $t_{c}=(2n+1 +2c)^{\frac{1}{2}} +  2^{-\frac{1}{2}}\,
n^{-\frac{1}{6}}\,s$,
\begin{equation}
F_{n+1,2}^{G}(t_{c})=F_{2}(s) \left\{ 1 + a_{c,2}^{G}\,u_{0}(s)\,
n^{-\frac{1}{3}} +b_{2}^{G}\,E_{c,2}^{G}(s)\,n^{-\frac{2}{3}}
\right\} + O(n^{-1})
\end{equation}
\end{lemma}
For $c=\frac{1}{2}$ we have faster convergence in Lemma
\ref{edgeworth for hermite} since $\; a_{\frac{1}{2},2}^{G}=0$.

\subsection{Edgeworth Expansion of $F_{n,2}^{L}$ }

From the similarities between the two expressions in
\eqref{hermitekernel} and \eqref{laguerrekernel}, the derivation
of $F_{n+1,2}^{L}$ follows exactly the same steps as the previous
one, they would differ only by some constant terms. \\

We see that the corresponding inner products in terms of the
$\alpha$'s are
\begin{eqnarray*}
& &
(\alpha_{1},\beta_{1})=\frac{2-c}{2^{\frac{1}{3}}}u_{0}n^{-\frac{1}{3}}\quad,
\quad
(\alpha_{5},\beta_{5})=-\frac{n^{-\frac{2}{3}}}{2^{\frac{2}{3}}\cdot
10}2w_{1} \quad , \quad
(\alpha_{6},\beta_{6})=-\frac{n^{-\frac{2}{3}}}{2^{\frac{2}{3}}\cdot
10}2w_{1},\\
& &
(\alpha_{2},\beta_{2})=\frac{n^{-\frac{2}{3}}}{2^{\frac{2}{3}}\cdot
10}2u_{2}\quad , \quad
(\alpha_{3},\beta_{3})=\frac{n^{-\frac{2}{3}}}{2^{\frac{2}{3}}\cdot
10}2(u_{2} -u_{1}v_{0} +u_{0}v_{1}-u_{0}v_{0}^{2}+u_{0}^{2}w_{0}),\\
& &
(\alpha_{4},\beta_{4})=\frac{n^{-\frac{2}{3}}}{2^{\frac{2}{3}}\cdot
10}2u_{2} \quad , \quad
(\alpha_{7},\beta_{7})=-\frac{n^{-\frac{2}{3}}}{2^{\frac{2}{3}}\cdot
10}(18-20c+5c^2)v_{0}= (\alpha_{8},\beta_{8}).
\end{eqnarray*}
If we set
\begin{equation*}
a_{c,2}^{L}(s)=\frac{c-2}{2^{\frac{1}{3}}} \;, \quad
b^{L}_{2}=\frac{2^{\frac{1}{3}}}{10} \quad \textrm{and}
\end{equation*}
\begin{equation}\label{second term}
E_{c,2}^{L}(s)=2w_{1}-3u_{2}+(5c^2 -20c +18)v_{0} +
u_{1}v_{0}-u_{0}v_{1} +u_{0}v_{0}^2   -u_{0}^{2}w_{0},
\end{equation}
our formula reads
\begin{lemma}\label{edgeworth for laguerre}
For $t_{c}= 4n +2\alpha +2c + 2(2n)^{\frac{1}{3}}s$,
\begin{equation}
F_{n+1,2}^{L}(t_{c})=F_{2}(s) \{ 1 + a_{c,2}^{L}\,u_{0}(s)\,
n^{-\frac{1}{3}} + b_{2}^{L}\,E_{c,2}^{L}(s)\,n^{-\frac{2}{3}} \}
\quad + \quad O(n^{-1}).
\end{equation}
\end{lemma}
For $c=2$ we obtain a faster convergence  as $a_{2,2}^{L}=0$.

\subsection{Fine tuning}\label{fineTuning}
To complete this analysis, we need to find values for the constant
$\;c=c_{1}\;$ in Lemma \ref{edgeworth for hermite} and
$\;c=c_{2}\;$ in Lemma \ref{edgeworth for laguerre} for which
$\quad E_{c_{1},2}^{G}(s)=E_{c_{2},2}^{L}(s)$. Which is equivalent
to
\begin{equation}\label{ellipse}
-20c_{1}^{2} +20c_{1}-2=5c_{2}^{2}-20c_{2}+18 \quad
\Leftrightarrow \quad (c_{1}-\frac{1}{2})^{2}
+\frac{1}{4}(c_{2}-2)^{2}=\frac{1}{4}
\end{equation}
This suggest the following change of variables.
\begin{equation*}
c_{G}=c_{1}-\frac{1}{2} \quad \textrm{and} \quad
c_{L}=\frac{c_{2}-2}{2}. \quad \textrm{Therefore}
\end{equation*}
\begin{equation*}
\eqref{ellipse} \quad \textrm{ changes to }\quad
c_{G}^{2}+c_{L}^{2}=\frac{1}{4}\;,\quad a_{c_{G},2}^{G}=c_{G} \, ,
\quad a_{c_{L},2}^{L}=2^{\frac{2}{3}}c_{L} \; ,
\end{equation*}
\begin{equation*}
-20c_{1}^{2}+20c_{1}-2 \;= \; -20c_{G}^{2}+3 \;, \quad
 5c_{2}^{2}-20c_{2}+18 \;=\; 20c_{L}^{2}-2 \; ,
\end{equation*}
\begin{equation}\label{finetune scale}
  t_{c_{G}}=(2(n+1)
+2c_{G})^{\frac{1}{2}}+2^{-\frac{1}{2}}n^{-\frac{1}{6}}s\quad
\textrm{and}\quad
 t_{c_{L}}=
4(n+1) +2\alpha + 4c_{L} +2(2n)^{\frac{1}{3}}s \; .
\end{equation}
We can now give a scaling of $t$ in terms of the size of the
matrices.
\begin{equation*}
  t_{c_{G}}=(2(n+1)
+2c_{G})^{\frac{1}{2}}+2^{-\frac{1}{2}}(n+1)^{-\frac{1}{6}}s\quad
\textrm{for the Gaussian case, and}
\end{equation*}
\begin{equation*}
 t_{c_{L}}=
4(n+1) +2\alpha + 4c_{L} +2(2(n+1))^{\frac{1}{3}}s \quad
\textrm{for the Laguerre case.}
\end{equation*}
Since all the functions derived so far are all differentiable,
this new scaling will only change the error function but not its
order and class. To keep the notations light, we will use the same
variable to represent these error functions. We therefore have
Theorem \ref{Hermite kernel}  and Theorem \ref{Laguerre kernel}
for an ensemble of $(n+1) \times (n+1)$ matrices.

\begin{figure}[h]\label{Comp1pther}
\includegraphics[height=75mm]{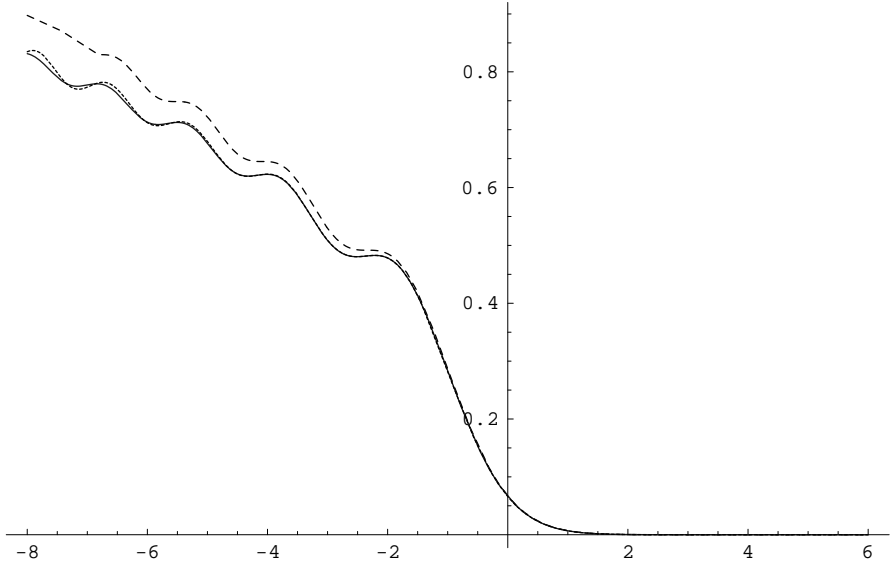}
\caption{For $c_{G}=0$ and $n=40$, The dashed curve is the usual
approximation of the one point correlation function from the Airy
kernel (which is the first two terms in \eqref{onepointher}), the
solid curve is the exact scaled one point correlation function and
the doted curve is our approximation \eqref{onepointher} }
\end{figure}

\begin{thm}\label{Gaussian kernel}
For $\; x  =  (2(n+1)
+2c_{G})^{\frac{1}{2}}+2^{-\frac{1}{2}}(n+1)^{-\frac{1}{6}}X\;$
and \\ $ \;y = (2(n+1)
+2c_{G})^{\frac{1}{2}}+2^{-\frac{1}{2}}(n+1)^{-\frac{1}{6}}Y$
\begin{eqnarray*}
K_{n+1}(x,y)\,dx =  \left\{  K_{\airy}(X,Y) -c_{G}\airy(X)\airy(Y)
(n+1)^{-\frac{1}{3}} +   \right.
\end{eqnarray*}
\begin{eqnarray*}
\frac{1}{20}\left[(X + Y)\airy^{'}(X)\airy ^{'}(Y) -
(X^2+XY+Y^2)\airy(X)\airy(Y) + \right.
\end{eqnarray*}
\begin{equation}\label{gaussian kernel}
\left. \left. \frac{-20c_{G}^2 +3 }{2}(\airy ^{'}(X)\airy(Y) +
\airy(X)\airy ^{'}(Y)) \right](n+1)^{-\frac{2}{3}}  +O((n+1)^{-1})
E(X,Y) \right\}  dX
\end{equation}
The error term, $ E(X,Y)$, is again a kernel of an integral
operator on $L^{2}(J)$ which is trace class for any subset $J$ of
the reals that is bounded away from minus infinity.
\end{thm}
Taking the limit as $Y\rightarrow X$ in \eqref{gaussian kernel}
give the one point correlation function $\rho_{n+1}$.
\begin{cor}
\begin{equation*}
\textrm{For} \quad x =(2(n+c_{G}))^{\frac{1}{2}} +
\frac{X}{2^{\frac{1}{2}}n^{\frac{1}{6}}},
\end{equation*}
\begin{equation*}
2^{-\frac{1}{2}}n^{-\frac{1}{6}}\rho_{n}(x)=2^{-\frac{1}{2}}n^{-\frac{1}{6}}K_{n}(x,x)=
[\airy^{'}(X)]^{2} -X [\airy(X)]^{2}
-c_{G}\,[\airy(X)]^{2}n^{-\frac{1}{3}} +
\end{equation*}
\begin{equation}\label{onepointher}
\frac{1}{20}\left\{  2X[\airy^{'}(X)]^{2} -3X^{2}[\airy(X)]^{2} +
(3-20c_{G}^{2})\airy^{'}(X)\airy(X) \right \} n^{-\frac{2}{3}} +
O(n^{-1})F_{n}(X)
\end{equation}
\end{cor}
Note that for $c_{G}=0$, this is formula (72) of \cite{Forr1}.
Figure \ref{Comp1pther} illustrates the accuracy of equation
\eqref{onepointher}.

\begin{figure}[h]
\includegraphics[height=75mm]{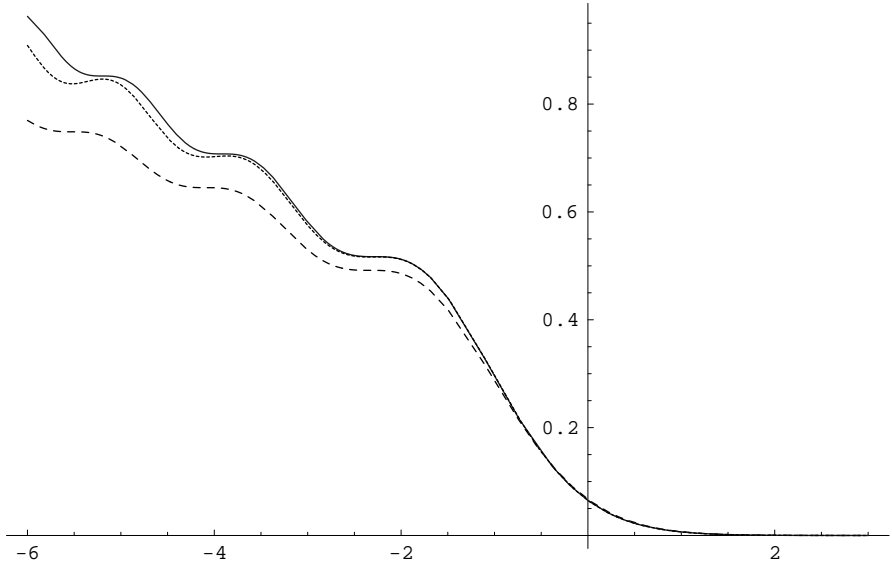}
\caption{For $c_{L}=0 ,\ \alpha =\frac{1}{2}$ and $n=40$, the
solid curve is the scaled one point correlation function, the
doted curve is our approximation from \eqref{onepointlag}, and the
dashed one is the first order approximation from the Airy
kernel.\label{1ptcorftlaggph}}
\end{figure}

For the Laguerre case,
\begin{thm} \label{Lag kernel}
For $\; x =4(n+1) +2\alpha + 4c_{L} +2(2(n+1))^{\frac{1}{3}}X
 \quad$ and \\ $y =4(n+1) +2\alpha + 4c_{L} +2(2(n+1))^{\frac{1}{3}}Y\;$  with
 $\;
 X,c_{L}$
 and $\; Y \;$
 bounded,
 \begin{eqnarray*}
 K_{n+1}^{\alpha}(x,y)\, dx= \left\{  K_{\airy}(X,Y)
-2^{\frac{2}{3}}c_{L}\airy(X)\airy(Y) (n+1)^{-\frac{1}{3}} +
\right.
\end{eqnarray*}
\begin{eqnarray*}
\frac{2^{\frac{1}{3}}}{ 10} \left[ \quad (X^2 +XY
+Y^2)\airy(X)\airy(Y) -(X+Y)\airy^{'}(X)\airy^{'}(Y) - \right.
\end{eqnarray*}
\begin{equation}\label{lag kernel}
\left. \left. (10c_{L}^2-1)(\airy(X)\airy^{'}(Y)
+\airy^{'}(X)\airy(Y)) \right](n+1)^{-\frac{2}{3}} +
O((n+1)^{-1})F(X,Y) \right\} \,dX
\end{equation}
The error term $F(X,Y)$ is again the kernel of an integral
operator on $L^{2}(J)$ which is trace class for any subset $J$ of
the reals which is bounded away from minus infinity.
\end{thm}
Taking the limit as $y$ goes to $x$ in \eqref{lag kernel} gives
the one point correlation function $\rho_{n+1}^{\alpha}$ in the
Laguerre case.
\begin{cor}
For $ x =4(n+c_{L})+2\alpha + 2(2n)^{\frac{1}{3}}X \, ,$
\begin{equation*}
2(2n)^{\frac{1}{3}}\rho_{n}^{\alpha}(x) = [\airy^{'}(X)]^{2}
-X[\airy(X)]^{2}- 2^{\frac{2}{3}}\, c_{L}[\airy(X)]^{2}
n^{-\frac{1}{3}} +
\end{equation*}
\begin{equation}\label{onepointlag}
\frac{2^{\frac{1}{3}}}{10}[  3X^2 [\airy(X)]^{2}
-2X[\airy^{'}(X)]^{2} + 2(1 -10\,c_{L}^2)\airy(X)\airy^{'}(X)
]n^{-\frac{2}{3}} + O(n^{-1})F(X)
\end{equation}
\end{cor}
Note that for $\; c_{L}=-\alpha /2 \;$ this is formula (73) in
\cite{Forr1}.
Figure \ref{1ptcorftlaggph} illustrates the accuracy of our
result. Combining Lemma \ref{edgeworth for hermite}, Lemma \ref{edgeworth
for laguerre} and the fine-tuned constants  in \eqref{finetune
scale} give Theorem \ref{Edgeworth expansion}.

\pagebreak

\appendix

\section{Proof of Theorem\ref{plancherel for Laguerre}}
 In this section we will follow Szeg\"o \cite{Szeg1}, Section 8.75.
 With the following changes.

\begin{itemize}
    \item We introduce a new variable $\;c\;$ in the definition of $\;l_{n}\;$
    to fine-tune the final result. $\;l_{n}=(4n +2\alpha
+2c)^{\frac{1}{2}} \;$instead of $\; l_{n}=(4n +2\alpha
+2)^{\frac{1}{2}}$.
    \item We will use the second order Hankel expansion of the Bessel
    function of large argument instead of the first order.
    \item We define $\; \xi = l_{n}
    +(2l_{n})^{-\frac{1}{3}}t \;$ instead of $\;\xi = l_{n}
    -(6l_{n})^{-\frac{1}{3}}t \;$
    to match the definition of the Airy function with the
    one commonly used.
    \item We give an estimate of the error term as a
    function of the independent variable
\end{itemize}
recalling the generating function of the Laguerre polynomials.
\begin{equation}\label{genfunct}
\sum_{n=0}^{\infty}\frac{L_{n}^{\alpha}(x)}{\Gamma(n+\alpha+1)}\omega^{n}=
e^{\omega}(x\omega)^{-\frac{\alpha}{2}}J_{\alpha}(2(x\omega)^{\frac{1}{2}})
\end{equation}
 The substitutions $x \rightarrow
\xi$ and $ \omega \rightarrow -\frac{\omega^{2}}{4}$, and for edge
scaling  $l_{n}=(4n +2\alpha +2c)^{\frac{1}{2}}$ ,
  $\omega \rightarrow l_{n}z .\;$  we deduce from \eqref{genfunct}

\begin{equation}\label{eq.1}
\frac{L_{n}^{\alpha}(\xi^{2})}{\Gamma(n+\alpha+1)}(-\frac{1}{4})^{n}=
\frac{2^{\alpha}(\xi)^{-\alpha}}{2 \pi i} l_{n}^{-2n -\alpha
}\int_{\gamma}e^{-\frac{l_{n}^{2}z^{2}}{4}
-\frac{1}{2}l_{n}^{2}log(z)}z^{c-1} e^{\frac{\pi \alpha
i}{2}}J_{\alpha}(e^{-\frac{\pi i}{2}}\xi l_{n}z) dz \quad
\end{equation}
where $\gamma$ is a symmetric contour enclosing the origin.

Using the Hankel expansion of the Bessel function of large
argument (see, for
example, Olver \cite[pgs 130--132]{Olve1} ): that\footnote{In the following
formula the $O$-term is actually of the form
\[O(\frac{1}{|\xi l_{n} z |^2})e^{\xi l_{n}z} + O(\frac{1}{|\xi
l_{n} z |^2})e^{-\xi l_{n}z + (\alpha + \frac{1}{2})\pi i}. \]} as
$|l_{n}z|\rightarrow \infty$ in $|\arg(z)| \leq \pi -\delta (< \pi
)$,

\begin{eqnarray*}
 e^{\frac{\pi \alpha i}{2}}J_{\alpha}(e^{-\frac{\pi i}{2}}\xi
l_{n}z)   &=&  (2 \pi \xi
             l_{n})^{-\frac{1}{2}}z^{-\frac{1}{2}}\left[e^{\xi l_{n}z} + e^{-\xi
              l_{n}z + (\alpha + \frac{1}{2})\pi i}  \right.\\
          &+& \left. \frac{4 \alpha^2 -1}{8 \xi l_{n} z}( -e^{\xi l_{n} z} +
             e^{-\xi l_{n} z + (\alpha + \frac{1}{2})\pi i}) +O(\frac{1}{|\xi
            l_{n} z |^2})\right]
\end{eqnarray*}

Substituting this in \eqref{eq.1} for $\xi = l_{n}
+(2l_{n})^{-\frac{1}{3}}t \; $  with $t$ bounded,

\begin{equation*}
\frac{L_{n}^{\alpha}(\xi^{2})}{\Gamma(n+\alpha+1)}(-\frac{1}{4})^{n}
  =
 \frac{2^{\alpha}(l_{n} +(2l_{n})^{-\frac{1}{3}}t)^{-\alpha}}{2 \pi i} l_{n}^{-2n -\alpha
}(2 \pi (l_{n} +(2l_{n})^{-\frac{1}{3}}t)l_{n})^{-\frac{1}{2}}
\cdot
\end{equation*}
\begin{equation*}
 \left\{\int_{\gamma}(1-\frac{4 \alpha^2 -1}{8
l_{n}^{2}
z(1+2^{-\frac{1}{3}}l_{n}^{-\frac{4}{3}}t)})e^{-\frac{l_{n}^{2}z^{2}}{4}
-\frac{1}{2}l_{n}^{2}\log(z) +l_{n}^{2}z  +2^{-\frac{1}{3}}l_{n}^{-
\frac{1}{3}}l_{n}z t} \quad z^{c-\frac{3}{2}}  dz   \quad+\right.
\end{equation*}
\begin{equation*}
   e^{(\alpha +\frac{1}{2})\pi i}
 \int_{\gamma}(1+\frac{4 \alpha^2
-1}{8 l_{n}^{2}
z(1+2^{-\frac{1}{3}}l_{n}^{-\frac{4}{3}}t)})e^{-\frac{l_{n}^{2}z^{2}}{4}
-\frac{1}{2}l_{n}^{2}\log(z) -l_{n}^{2}z -2^{-\frac{1}{3}}l_{n}^{-
\frac{1}{3}}l_{n}zt} \quad z^{c-\frac{3}{2}}  dz
\end{equation*}
\begin{equation}\label{eq.2}
  + \left.
0(l_{n}^{-4})\int_{\gamma}|e^{-\frac{l_{n}^{2}z^{2}}{4}
-\frac{1}{2}l_{n}^{2}\log(z) \pm l_{n}^{2}z \pm
2^{-\frac{1}{3}}l_{n}^{- \frac{1}{3}}l_{n}zt}  | \quad |dz | \quad
\quad \right\}
\end{equation}

The $\pm$ signs would be taken according to which one gives the
larger contribution.  To simplify notations we set $p= (4
\alpha^2 -1)/(8 l_{n}^{2}
z(1+2^{-\frac{1}{3}}l_{n}^{-\frac{4}{3}}t))$.\\

The integral over $\gamma$ can be split into the path in the upper
half-plane and the path in the lower half-plane. However the
lower-half plane part of the path can be transformed into the
upper half plane via the transformation $z \rightarrow \bar{z}$,
taking into account the orientation change and the fact that the
integrand is an analytic function of $z$ on $\gamma$. The
integrand will be transformed into its complex conjugate under
this transformation. The Jacobian of the transformation is $-1$ so
it will reverse the orientation once again. Thus the lower half
plane contribution of the integral is equal to the conjugate of
the upper half contribution such that the integral over $\gamma$
is exactly twice the real part of the integral over the upper half
portion $ \gamma_{+}$ of the contour. Thus \eqref{eq.2} is equal
to
\begin{equation}\label{laguerre}
2^{\alpha}(l_{n} +(2l_{n})^{-\frac{1}{3}}t)^{-\alpha} l_{n}^{-2n
-\alpha }(2 \pi (l_{n}
+(2l_{n})^{-\frac{1}{3}}t)l_{n})^{-\frac{1}{2}} \cdot
2\mathcal{R}e \{ \frac{1}{2 \pi i}G + \frac{1}{2 \pi i}H + K \}
\end{equation}
where
\begin{equation}\label{g}
 G =\int_{\gamma_{+}}(1-p)e^{-l_{n}^{2}f_{1}(z)
+ 2^{-\frac{1}{3}}l_{n}^{ \frac{2}{3}}z\,t} \, z^{c-\frac{3}{2}}\,
dz = \int_{\gamma_{+}}e^{-l_{n}^{2}f_{1}(z)}g(z)\,dz
\end{equation}
\begin{equation*}
\textrm{with}\quad f_{1}(z)=\frac{z^{2}}{4} -z + \frac{1}{2}\log(z)
\end{equation*}
\begin{equation}\label{h}
 H =e^{(\alpha + \frac{1}{2}) \pi i }\int_{\gamma_{+}}
(1+p)e^{-l_{n}^{2}f_{2}(z) -2^{-\frac{1}{3}}l_{n}^{
\frac{2}{3}}zt} \, z^{c-\frac{3}{2}}\, dz =\int_{\gamma_{+}}
e^{-l_{n}^{2}f_{2}(z)}h(z)\,dz
\end{equation}
\begin{equation*}
\textrm{with} \quad f_{2}(z)= \frac{z^{2}}{4}+z +
\frac{1}{2}\log(z)
\end{equation*}

\begin{equation}\label{k}
\textrm{and} \quad  K
=0(l_{n}^{-4})\int_{\gamma_{+}}|e^{-\frac{l_{n}^{2}z^{2}}{4}
-\frac{1}{2}l_{n}^{2}log(z) \pm l_{n}^{2}z \pm
2^{-\frac{1}{3}}l_{n}^{ \frac{2}{3}}zt}  | \quad |dz |
\end{equation}

We will use the steepest descent method to find an asymptotic
expansion for $G,H$ and $K$ for large $l_{n}$ keeping $\alpha$
fixe and $c$ bounded.\footnote{Actually $c$ can grow with $n$, but we
are not going to look at this problem.}

 \subsection{Steepest descent method for G}

The steepest descent for this section is with respect to large
$l_{n}$ in \eqref{g}. The saddle point condition $f^{'}(z)=0$
gives $z_{0}=1$ with $f^{''}(1)=0$.

The steepest descent curve leaves $z_{0} =1$ at angles $ 0 ,\
\pi/3,
 \ $ or $\  2\pi/3 $. The direction of maximum decreases of
 $f_{1}$ is $ \  2\pi/3$.  We deform
 $\gamma_{+}$  at $z_{0}=1$
such that the resulting contour leaves $z_{0}=1$ at angle
$2\pi/3\;$ as the line segment
\begin{equation}\label{eq.3}
z=1+2^{\frac{1}{3}}l_{n}^{-\frac{2}{3}} \rho e^{\frac{2 \pi i}{3}}
\quad \textrm{ with }   0 \leq \rho \leq n^{\delta} \quad \textrm{
where } \quad 0 < \delta < \frac{1}{6}
\end{equation}
then along the segment symmetric to this segment with respect to
the imaginary axis and finally connect the tip of these two line
segments by an arc of circle centered at the origin with radius
$r$. For simplicity, we will also call this path\footnote{See
\cite[Section 8.75]{Szeg1} for an illustration of this path.} $
\gamma_{+}$.

In the following sections, we will estimate the contribution of
each portion of the contour to $G$.

\subsubsection{On the arc of circle $ z= re^{\phi i}$\\\\}
 Since $\mathcal{R}e(-f_{1}(z))=-\frac{1}{4}r^{2} \cos2\phi
+r\cos\phi -\frac{1}{2}\log(r)$ is a decreasing function of $\phi$
for $0 < \phi < \pi$, the major contribution of this arc is
bounded above by the value of the integrand at the end point where
 $ \rho = n^{\delta}$.
 A deformation of the path of integration near $-1$ is immaterial
 so the estimate of the remainder of the path is again bounded
 above by the value of the integrand at $\rho =  n^{\delta}$.
 (Therefore in the next section we will focus on the asymptotics on
 the line segment in the first quadrant only.)

 If we set $\;z=re^{i \theta}\;$ where
 $\;r^{2}= 1+2^{\frac{2}{3}}l_{n}^{-\frac{4}{3}}n^{2 \delta} -
  2^{\frac{1}{3}}l_{n}^{-\frac{2}{3}}n^{\delta}\;$  and $\;
 \theta_{0} \leq \theta \leq \pi - \theta_{0} \;$  where $\;\theta_{0} \;$ is the
 angle that the ray from the origin to the tip of the line segment
 makes with the real axis,
 $\; \cos(\theta_{0})= \frac{ 2- 2^{\frac{1}{3}}l_{n}^{-\frac{2}{3}}n^{\delta}}{2r}.$
The modulus of this integrand in \eqref{g} is of order
\begin{equation*}
   r^{c-\frac{3}{2}}e^{k(\theta)} \quad  \textrm{where} \quad  k(\theta)=
 -\frac{1}{4}r^{2}\,\cos(2 \theta) + r\,\cos( \theta)
 -\frac{1}{2} \log(r)
 +2^{-\frac{1}{3}}l_{n}^{-\frac{4}{3}}\,r\,t \cos( \theta)
 \end{equation*}
$k^{'}$ has roots $\; \pi,\; 0 \;$  and
 $ \; \theta_{2} \; \textrm{where}\; \cos(\theta_{2})=
 \frac{1+2^{-\frac{1}{3}}l_{n}^{-\frac{4}{3}}t}{r} \;$,
 for large $\;l_{n}$, $\; 0 \leq \theta_{2} \leq  \theta_{0}\; $.
$k$ increases from $ 0$ to $\theta_{2}$ and decreases from
 $\theta_{2}$ to $\pi$.  Therefore the maximum of $k(\theta)$ on the arc of circle is
 at $\theta_{0}$. We can also use this as an upper bound of the
 contribution of the contour around $-1$.  Thus the contribution of the arc
 of circle and the line segment on the second quadrant is of order of the
 modulus of the integrand evaluated at the tip of the line
  segment in the first quadrant.We can therefore focus our attention in the next section
  on just the line segment in the first quadrant.  The
  parametrization \eqref{eq.3} shows that $\rho$ will go to
  infinity with $n$, so we want to estimate the contribution not only  of the line
  segment but of the whole ray $\;0 \leq \rho \leq \infty \;$. The
  error that we make by taking the ray is
 $\;| \int_{n^{\delta}}^{\infty} z^{c-\frac{2}{3}}
 e^{-l_{n}^{2}f_{1}(z)+ 2^{-\frac{1}{3}}l_{n}^{\frac{2}{3}}tz}
 dz |$  .\\
 On the ray, $ \;-l_{n}^{2}f_{1}(z)= \frac{3}{4} l_{n}^{2}
 -\frac{\rho^{3}}{3} +O(\rho^{4})l_{n}^{-\frac{2}{3}}\;$, thus\\
 $ \int_{n^{\delta}}^{\infty} z^{c-\frac{2}{3}}
 e^{-l_{n}^{2}f_{1}(z)+ 2^{-\frac{1}{3}}l_{n}^{\frac{2}{3}}tz}
 dz =O(1) 2^{\frac{1}{3}}l_{n}^{-\frac{2}{3}} e^{\frac{3}{4}
 l_{n}^{2} + 2^{-\frac{1}{3}}l_{n}^{\frac{2}{3}}t}
 \int_{n^{\delta}}^{\infty}e^{-\frac{\rho^3}{3} + \rho t
 e^{\frac{2 \pi i}{3}}} d\rho$ .\\
We therefore need to give an estimate an estimate of the integral
on the right of this last equality and an estimate of the
integrand when $\;\rho= n^{\delta}$ in \eqref{g}.\\

\underline{Integrand when $\;\rho= n^{\delta}$}\\
The order of the integrand for large $n$ is $\quad O(1)
2^{\frac{1}{3}}l_{n}^{-\frac{2}{3}} e^{\frac{3}{4}
 l_{n}^{2} + 2^{-\frac{1}{3}}l_{n}^{\frac{2}{3}}t}
 \, e^{-\frac{\rho^3}{3} + \rho t
 e^{\frac{2 \pi i}{3}}}.$\\
We need to give an estimate of the factor $\;u:=
e^{-\frac{\rho^3}{3} + \rho t
 e^{\frac{2 \pi i}{3}}}\;$ as a function of $\;t$.\\
 If $t$ is positive, it is of order $e^{-\frac{n^{3 \delta}}{3}}
 \cdot e^{-\frac{t n^{\delta}}{2}}$ or of order $e^{-\frac{n^{3 \delta}}{3}}
 \cdot e^{-t}.$\\
 If $t$ is negative, we have
$\; |u|=e^{-\frac{\rho^3}{3} + \frac{\rho y}{2}}$ where $y=-t$.
For trace class convergence we need an estimate that will decay
exponentially for large $t$. The expansion of the exponent in
$|u|$ around his critical point $\rho_{0}= - \sqrt{\frac{y}{2}}$
is
 $-\frac{2}{3}(\frac{y}{2})^{\frac{3}{2}} +\sqrt{\frac{y}{2}}(\rho
 +\sqrt{\frac{y}{2}})^2 -\frac{1}{3}(\rho +\sqrt{\frac{y}{2}})^3$.
Thus for $\rho =n^{\delta}, \;  u$ is of order
 $e^{-\frac{(n^{ \delta}+ \sqrt{\frac{y}{2}})^3}{3}}
 \cdot e^{-\frac{2}{3}(\frac{y}{2})^{\frac{3}{2}}}$ or of order
 $e^{-\frac{n^{3 \delta}}{3}}\cdot e^{-y}$ since $y$ is bounded and
 positive.
 So in either case the contribution is of order,
 $e^{-\frac{n^{3 \delta}}{3}}\cdot e^{-|t|}$.\\\\

\underline{Tail integration}\\
 We show that the contribution of
 $\int_{n^{\delta}}^{\infty}e^{-\frac{\rho^3}{3} -\frac{\rho t}{2}
 } d\rho \quad$ is of the same order.\\

 If $\; t \geq 0\, ,$
$\quad \int_{n^{\delta}}^{\infty}e^{-\frac{\rho^3}{3} -\frac{\rho
t}{2}
 } d\rho \, \leq \int_{n^{\delta}}^{\infty}e^{-\frac{\rho^3}{3} - \frac{t}{2}
 } d\rho \,  = e^{-\frac{t}{2}}\int_{n^{\delta}}^{\infty}e^{-\frac{\rho^3}{3}
 } d\rho \, \leq  e^{-\frac{t}{2}}\int_{n^{\delta}}^{\infty}e^{-\frac{\rho}{3}
 } d\rho $\\
so it is of order $e^{-\frac{n^{\delta}}{3}}\cdot
e^{-\frac{t}{2}}$\\
If $t$ is negative, a similar change of variable and expansion
of the integrand leads to\\
$\int_{n^{\delta}}^{\infty}e^{-\frac{\rho^3}{3} -\frac{\rho t}{2}
 } d\rho \quad = \int_{n^{\delta}}^{\infty}e^{-\frac{2}{3}(\frac{y}{2})^{\frac{3}{2}} +\sqrt{\frac{y}{2}}(\rho
 +\sqrt{\frac{y}{2}})^2 -\frac{1}{3}(\rho +\sqrt{\frac{y}{2}})^3}
 d\rho $\\
 $e^{-\frac{2}{3}(\frac{y}{2})^{\frac{3}{2}}}
 \int_{n^{\delta}}^{\infty}e^{  -\frac{1}{3}(\rho
+\sqrt{\frac{y}{2}})^3[1-3(\rho+\sqrt{\frac{y}{2}})^{-1}\sqrt{\frac{y}{2}}]}
 d\rho \quad \sim e^{-\frac{2}{3}(\frac{y}{2})^{\frac{3}{2}}}
 \int_{n^{\delta}}^{\infty}e^{  -\frac{1}{3}(\rho
+\sqrt{\frac{y}{2}})^3}
 d\rho$ \\
 This is of order, $e^{-\frac{2}{3}(\frac{y}{2})^{\frac{3}{2}}}
 \int_{n^{\delta}}^{\infty}e^{  -\frac{1}{3}(\rho
+\sqrt{\frac{y}{2}})}
 d\rho$ or of order $ e^{-\frac{2}{3}(\frac{y}{2})^{\frac{3}{2}}}
 \cdot e^{  -\frac{1}{3}(n^{\delta}
+\sqrt{\frac{y}{2}})}$\\
So the contribution of this integral is of order $e^{-\frac{n^{
\delta}}{3}}\cdot e^{-|t|}$

 We conclude this subsection by recording that the error that we
 make by neglecting the remainder of the contour and considering
 the integral from zero to infinity instead of zero to
 $n^{\delta}$ is at most of order

 \begin{equation}\label{remainder}
 2^{\frac{1}{3}}l_{n}^{-\frac{2}{3}}e^{\frac{3}{4}l_{n}^{2}
 +2^{-\frac{1}{3}}l_{n}^{\frac{2}{3}}t } e^{-\frac{ n^{
 \delta}}{3}}\cdot e^{-|t|}
 \end{equation}\\

\subsubsection{On the ray $z= 1 +
 2^{\frac{1}{3}}l_{n}^{-\frac{2}{3}} \rho e^{\frac{2 \pi
 i}{3}}, \quad \rho \in[0,\infty)$}
The Taylor expansion of $f_{1}$ at $z_{0}=1$ is\\
 $ \quad f_{1}(z)=-\frac{3}{4} +\frac{ \rho^{3}}{3} l_{n}^{-2} -
\frac{1}{2}
 \sum_{k=4}^{\infty} c_{1,k}( \rho l_{n}^{-\frac{2}{3}})^{k} \quad$ with
 $\; c_{1,k}=(-1)^{k} \frac{2^{\frac{k}{3}}e^{\frac{2 \pi
 k i}{3}}}{k}$ ,  and\\

   $\quad 2^{-\frac{1}{3}}l_{n}^{\frac{2}{3}}tz =
 2^{-\frac{1}{3}}l_{n}^{\frac{2}{3}}t + \rho t e^{\frac{2 \pi i
 }{2}}.\quad$  Taking in account the error estimate from \eqref{remainder}, the substitution of
 these quantities in $G$ give,

\begin{equation*}
 G =
 2^{\frac{1}{3}}l_{n}^{-\frac{2}{3}}e^{\frac{2 \pi i}{3}}e^{\frac{3}{4}l_{n}^{2}
 +2^{-\frac{1}{3}}l_{n}^{\frac{2}{3}}t }\cdot \left\{ \int_{0}^{\infty}
 e^{-\frac{\rho^3}{3}+\rho t e^{\frac{2 \pi
 i}{3}}} g_{n}(\rho) \, d\rho  \quad  +  O(e^{-\frac{ n^{
 \delta}}{3}})\cdot e^{-|t|} \right\} \quad \textrm{ with}
 \end{equation*}
 \begin{equation}\label{eq.4}
 g_{n}(\rho)=(1-\frac{4 \alpha^2 -1}{8 l_{n}^{2}
( 1 + 2^{\frac{1}{3}}l_{n}^{-\frac{2}{3}} \rho e^{\frac{2 \pi
i}{3}}
 )(1+2^{-\frac{1}{3}}l_{n}^{-\frac{4}{3}}t)})( 1 +
 2^{\frac{1}{3}}l_{n}^{-\frac{2}{3}} \rho e^{\frac{2 \pi i}{3}}
)^{c-\frac{3}{2}} e^{\frac{1}{2}
  \sum_{k=4}^{\infty}c_{1,k}(\rho l_{n}^{-\frac{2}{3}})^{k}l_{n}^{2}}
\end{equation}

 \subsection{Steepest descent for H}

 The analysis for $H$ differs from that of $G$ in the location of
 the saddle point, and the orientation of the contour. A similar
 analysis shows that the saddle point is now at $z_{0}=-1$, the
 final contour of integration is the same but oriented in the opposite
 direction. It leaves $z_{0}$ at angle $\pi/3$. The error estimate
 on the arc of circle and on the tail of the corresponding ray is the
 same. The new parametrization on the ray is $ \; z=-1 +
 2^{\frac{1}{3}}l_{n}^{-\frac{2}{3}}e^{\frac{\pi i}{3}} \rho = -(1 +
 2^{\frac{1}{3}}l_{n}^{-\frac{2}{3}}e^{-\frac{2 \pi i}{3}} \rho ),\quad
  0\leq\rho\leq \infty $\\
The Taylor expansion of $f_{2}$ at $z=-1$ is
 \[
f_{2}(z)=-\frac{3}{4} + \frac{\pi}{2}i +l_{n}^{-2}
\frac{\rho^{3}}{3} -
 \frac{1}{2} \sum_{k \geq 4} c_{2,k}
 \rho^{k}l_{n}^{-\frac{2k}{3}},\quad \textrm{with}\quad
 c_{2,k}=(-1)^{k}\frac{2^{\frac{k}{3}}e^{-\frac{2 \pi k i}{3}}}{k}
\]
This leads to
\begin{equation*}
 H =
 -2^{\frac{1}{3}}l_{n}^{-\frac{2}{3}}e^{-\frac{2 \pi i}{3}}e^{\frac{3}{4}l_{n}^{2}
 +2^{-\frac{1}{3}}l_{n}^{\frac{2}{3}}t }\cdot \left\{ \int_{0}^{\infty}
 e^{-\frac{\rho^3}{3}+\rho t e^{-\frac{2 \pi
 i}{3}}} h_{n}(\rho) \, d\rho  \quad  +  O(e^{-\frac{ n^{
 \delta}}{3}})\cdot e^{-|t|} \right\} \quad \textrm{ with}
 \end{equation*}
 \begin{equation}\label{eq.5}
 h_{n}(\rho)=(1-\frac{4 \alpha^2 -1}{8 l_{n}^{2}
( 1 + 2^{\frac{1}{3}}l_{n}^{-\frac{2}{3}} \rho e^{-\frac{2 \pi
i}{3}}
 )(1+2^{-\frac{1}{3}}l_{n}^{-\frac{4}{3}}t)})( 1 +
 2^{\frac{1}{3}}l_{n}^{-\frac{2}{3}} \rho e^{-\frac{2 \pi i}{3}}
)^{c-\frac{3}{2}} e^{\frac{1}{2}
  \sum_{k=4}^{\infty}c_{2,k}(\rho l_{n}^{-\frac{2}{3}})^{k}l_{n}^{2}}
\end{equation}

\subsection{Asymptotics for K}

The asymptotics of the integral factor in $\textbf{K}$ depends on
the leading term in the expansion of either $\textbf{G}$ or
$\textbf{H}$, depending on which one is larger as shown in
\eqref{k}. But from the previous analysis, the leading term of
both $\textbf{G}$ and $\textbf{H}$ are of the same order. Thus
$\textbf{K}$ is also of order of a nonzero linear combination of
the leading terms in $\textbf{G}$ and in $\textbf{H}$ times
$O(l_{n}^{-4})$. In our case\footnote{The choice of this
representation of the error is for trace class convergence of the
final result.} we take the linear combination to be $\quad
\mathcal{R}e\frac{1}{2\pi i} (\textbf{G} + \textbf{H} )$.

 \subsection{Conclusion}

 Note that from \eqref{eq.4}  and \eqref{eq.5} we see
 that except for the $O$-term, $H=-\overline{G}$.\\
The change of variable  $\;\rho \mapsto \rho e^{-\frac{2 \pi
i}{3}}\;$ transform
 $g_{n}(\rho)$ into a real function $g_{1}(\rho)$  and

\begin{eqnarray*}
\textbf{G}=
2^{\frac{1}{3}}l_{n}^{-\frac{2}{3}}e^{\frac{3}{4}l_{n}^{2}
 +2^{-\frac{1}{3}}l_{n}^{\frac{2}{3}}t }\cdot
 (\int_{0}^{\infty e^{\frac{2 \pi i}{3}}}e^{(-\frac{\rho^3}{3}+\rho t
  )}  g_{1}(\rho)d\rho + O(e^{-\frac{n^{\delta}}{3}}) e^{-|t|} ),
  \quad \textrm{note also that}
\end{eqnarray*}
\begin{equation*}
2\mathcal{R}e \frac{1}{2 \pi i}(G+H) =
2\cdot2^{\frac{1}{3}}l_{n}^{-\frac{2}{3}}e^{\frac{3}{4}l_{n}^{2}
 +2^{-\frac{1}{3}}l_{n}^{\frac{2}{3}}t }\cdot
\end{equation*}
\begin{equation}\label{eq.6}
\left\{ \frac{1}{2 \pi i}  [\int_{\infty e^{-\frac{2 \pi
i}{3}}}^{0}e^{(-\frac{\rho^3}{3}+\rho t )}  g_{1}(\rho)d\rho  +
\int_{0}^{\infty e^{\frac{2 \pi i}{3}}}e^{(-\frac{\rho^3}{3}+\rho
t )}  g_{1}(\rho)d\rho ] +O(e^{-\frac{n^{\delta}}{3}}) e^{-|t|}
\right\}
\end{equation}
At this point we can give an estimate for $\textbf{K}$ based on
this last formula since the leading term of $g_{1}$ is 1.
\begin{equation*}
\textbf{K}=O(l_{n}^4)\cdot
2^{\frac{1}{3}}l_{n}^{-\frac{2}{3}}e^{\frac{3}{4}l_{n}^{2}
 +2^{-\frac{1}{3}}l_{n}^{\frac{2}{3}}t }\cdot
 \frac{1}{2 \pi i}(\int_{\infty e^{-\frac{2 \pi i}{3}}}^{0}e^{(-\frac{\rho^3}{3}+\rho t
  )}  d\rho +\int_{0}^{\infty e^{\frac{2 \pi i}{3}}}e^{(-\frac{\rho^3}{3}+\rho t
  )}  d\rho)
 \end{equation*}
\begin{equation}\label{ErrorK}
= \quad O(l_{n}^{-4})\cdot
2^{\frac{1}{3}}l_{n}^{-\frac{2}{3}}e^{\frac{3}{4}l_{n}^{2}
 +2^{-\frac{1}{3}}l_{n}^{\frac{2}{3}}t }\cdot \airy(t)
\end{equation}
\begin{equation}\label{Airy}
\textrm{where} \quad  \airy(t)= \frac{1}{2 \pi i}(\int_{\infty
e^{-\frac{2 \pi i}{3}}}^{0}e^{(-\frac{\rho^3}{3}+\rho t
  )}  d\rho +\int_{0}^{\infty e^{\frac{2 \pi i}{3}}}e^{(-\frac{\rho^3}{3}+\rho t
  )}  d\rho)\quad \textrm{is the Airy function}
\end{equation}

To simplify notation, we will combine the two paths of integration
in $\textbf{G}$ and call the new path $\sigma$.

 With this notation, \eqref{laguerre} gives
\begin{equation*}
e^{-\frac{\xi^{2}}{2}}L_{n}^{\alpha}(\xi^{2})=  (-1)^{n} 2^{2n}
\Gamma(n+\alpha+1)
 2^{\alpha}(\xi)^{-\alpha}e^{-\frac{\xi^{2}}{2}} l_{n}^{-2n -\alpha
} (2 \pi \xi l_{n})^{-\frac{1}{2}}
2\cdot2^{\frac{1}{3}}l_{n}^{-\frac{2}{3}}e^{\frac{3}{4}l_{n}^{2}
 +2^{-\frac{1}{3}}l_{n}^{\frac{2}{3}}t } \cdot
\end{equation*}
\begin{equation}\label{Laguer}
 [\frac{1}{2 \pi i}\int_{\sigma}e^{(-\frac{\rho^3}{3}+\rho t
  )}  g_{1}(\rho)d\rho \quad +\quad O(l_{n}^{-4})Ai(t) \quad +
  \quad O(e^{-\frac{n^{\delta}}{3}}) e^{-|t|} ]
\end{equation}

Using \eqref{eq.4} and the help of {\em Mathematica}, we have the
following expansion\footnote{The expansion is valid for this
derivation up to $l_{n}^{-6\frac{2}{3}}$  as the error estimate
for $K$ indicates. We choose to stop here at
$l_{n}^{-4\frac{2}{3}}$ note also that the last term will served
to estimate the error.} of $g_{1}$ in powers of
$l_{n}^{-\frac{2}{3}}$.

\begin{equation*}
 1\; + \;
 2^{\frac{1}{3}}(\frac{{\rho}^{4}}{4}+ (c-\frac{3}{2}
 )\rho )\;l_{n}^{-\frac{2}{3}} \; +
2^{\frac{2}{3}} \left(-\frac{{\rho}^{5}}{5}+ \frac{{\rho}^{8}}{32}
+\frac{ (c -\frac{3}{2} ){\rho}^{5}}{4}
 +
 \frac{  (c-\frac{3}{2}
){\rho}^{2}(c-\frac{5}{2} )}{2} \right) \;l_{n}^{-\frac{4}{3}}
\end{equation*}
 \begin{equation*}
+\left(\frac{ (c-\frac{3}{2} ){\rho}^{6}  (c-\frac{5}{2} )}{4}
+\frac{\rho^{6}}{3} -\frac{\rho^{9}}{10}+  \frac {{ \rho}^{12}
}{192}+\frac{1}{3}(c-\frac{3}{2} ){\rho}^{3}(c-\frac{5}{2}
)(c-\frac{7}{2}) \quad - \right.
\end{equation*}
 \begin{equation*}
 \left.  \frac{a^{2}}{2}+\frac{1}{8} + (c-\frac{3}{2})2\rho(-\frac{{\rho}^{5}}{5}+\frac{{\rho}^{8}}{32}
 ) \right)l_{n}^{-2}\quad +
\end{equation*}
 \begin{equation*}
 \left(2^{\frac{1}{3}}(-\frac{2{\rho}^{7}}{7}+\frac{{\rho}^{10}}{12}+\frac{
 {\rho}^{10}}{25} -
\frac{{\rho}^{13}}{80}+\frac{{\rho}^{16}}{3072}) +\frac{1}{4}
\left(\frac{ (c-\frac{3}{2} ){\rho}^{3} (c-\frac{5}{2} )
(c-\frac{7}{2} )}{3} - \frac{a^{2}}{2} + \frac{1}{8} \right)\sqrt
[3]{2}{\rho}^{4} \right.
\end{equation*}
 \begin{equation*}
 + (c-\frac{3}{2} ){2}^{\frac{1}{3}}{\rho}^{2}
 (c-\frac{5}{2}
) (-\frac{{\rho}^{5}}{5}+\frac{{\rho}^{8}}{32} )
 + (c-\frac{3}{2} )2^{\frac{1}{3}}\rho
(\frac{{\rho}^{6}}{3}-\frac{\rho^{9}}{10}+\frac{\rho^{12}}{192}
 )
 \end{equation*}
 \begin{equation*}
 \left. +\frac{2^{\frac{1}{3}}}{12} (c-\frac{3}{2} ){\rho}^{4}  (c-\frac{5}{2} )
(c-\frac{7}{2} ) (c -\frac{9}{2} )
  + (-\frac{a^{2}}{2}+\frac{1}{8} )
(c-\frac{3}{2} )2^{\frac{1}{3}}\rho -\frac{1}{8} (-4\,{a}^{2}+1
)2^{\frac{1}{3}}\rho \right)l_{n}^{-8/3}
\end{equation*}
\begin{equation*}
+q(\rho)O(l_{n}^{-\frac{10}{3}}) \quad \textrm{for some
polynomial} \; q
\end{equation*}
 The integral in \eqref{Laguer} can be
expresses as a linear combination of the Airy function and its
derivative using $\;\airy^{(k)}(t)= \frac{1}{2 \pi
i}\int_{\sigma}e^{-\frac{\rho^3}{3}+\rho t
  }  \rho^k d\rho \;$ and this expansion of $\;g_{1}$. Using the
  Airy differential equation $\airy''(t)=t\airy(t)$, it reduces
  to an expression involving only the independent variable $t$, $\airy$ and $\airy'$ .\\
The contribution of the last term of the expansion of $g_{1}$ is a
finite combination of the form $\; \sum_{m=0}^{k}(p_{m}(t)
\,\airy(t) + q_{m}(t) \, \airy'(t)) \;$ for some polynomials $p$
and $q$. If $t$ is bounded away from minus infinity, this is of
order $\airy(t)$. In this paper we assume therefore that this is
the case for $t$.\\
If $\; \xi = l_{n} +(2l_{n})^{-\frac{1}{3}}t \;$ we have\\
$e^{-\frac{\xi^2}{2}}\xi^{-(\alpha
+\frac{1}{2})}=l_{n}^{-(\alpha+\frac{1}{2})}e^{-\frac{l_{n}^{2}}{2}
 -2^{-\frac{1}{3}}l_{n}^{\frac{2}{3}}t
 }[1 - \frac{t^2}{2^2}2^{\frac{1}{3}} l_{n}^{-\frac{2}{3}} +(
 \frac{t^4}{2^5}- (2\alpha+1)\frac{t}{2^2})2^{\frac{2}{3}}
 l_{n}^{-\frac{4}{3}}+ (
 -\frac{t^6}{3 \cdot 2^6}+ (2\alpha+1)\frac{t^3}{2^3})
 l_{n}^{-2}+ O(l_{n}^{-\frac{8}{3}} t^2)]$. \\
Stirling formula gives
 $\; \Gamma(n+\alpha +1)=2^{\frac{1}{2}}
 \pi^{\frac{1}{2}}n^{n+\alpha+\frac{1}{2}}e^{-n} [1+ \frac{6\alpha
 (\alpha +1) +1}{12n} +O(n^{-2})]\quad $ and  we have
$l_{n}^{-2n-2\alpha -\frac{5}{3}} = 2^{-2n-2\alpha
 -\frac{5}{3}}n^{-n-\alpha -\frac{5}{6}}e^{-\frac{\alpha +c}{2}}[1
 +\frac{(\alpha+c)(-9\alpha+3c-10)}{24n} +O(n^{-2})]$. \\\\
Thus\\
 $(-1)^{n} 2^{2n} \Gamma(n+\alpha+1)
 2^{\alpha}(\xi)^{-\alpha} l_{n}^{-2n -\alpha
} (2 \pi \xi l_{n})^{-\frac{1}{2}}
2\cdot2^{\frac{1}{3}}l_{n}^{-\frac{2}{3}}e^{\frac{3}{4}l_{n}^{2}
 +2^{-\frac{1}{3}}l_{n}^{\frac{2}{3}}t }e^{-\frac{\xi^2}{2}}= \\
 \\
 (-1)^{n}2^{-\alpha-\frac{1}{3}}n^{-\frac{1}{3}}(1 +
 \frac{3\alpha^2 +2\alpha -6\alpha c+3c^2 -10c +2}{24n}
 +O(n^{-2}))[1 - \frac{t^2}{2^2}2^{\frac{1}{3}} l_{n}^{-\frac{2}{3}} +(
 \frac{t^4}{2^5}- (2\alpha+1)\frac{t}{2^2})2^{\frac{2}{3}}
 l_{n}^{-\frac{4}{3}}+ (
 -\frac{t^6}{3 \cdot 2^6}+ (2\alpha+1)\frac{t^3}{2^3})
 l_{n}^{-2}+ O(l_{n}^{-\frac{8}{3}} t^2)]$.\\ \\
This in \eqref{Laguer} together with the expansion of $g_{1}$ give
the desired result.\footnote{Actually  \[\xi =(4n +2\alpha
+2c)^{\frac{1}{2}} + \frac{t}{2^{\frac{2}{3}}n^{\frac{1}{6}}} +
O(n^{-\frac{7}{6}}),\]  but due to the smoothness of
$\;e^{-\frac{\xi^2}{2}} L_{n}^{\alpha}(\xi^{2})$   the error
that we make by removing the $O$-term is negligible if we are
aiming for an accuracy of order $n^{-1}$.}

For $\xi =(4n +2\alpha +2c)^{\frac{1}{2}} +
\frac{t}{2^{\frac{2}{3}}n^{\frac{1}{6}}} \quad$ and $t\quad$
bounded,
\begin{eqnarray*}
 e^{-\frac{\xi^2}{2}}
L_{n}^{\alpha}(\xi^{2})=(-1)^{n}2^{-\alpha-\frac{1}{3}}n^{-\frac{1}{3}}
\{ \airy(t) +
 \frac{(c-1)}{2^{\frac{1}{3}}} \airy'(t) n^{-\frac{1}{3}} +
 \end{eqnarray*}
 \begin{eqnarray*}
[ \frac{2-10c+5c^2 -5 \alpha}{10 \cdot 2^{\frac{2}{3}} }t \airy(t)
+ \frac{t^2}{20 \cdot 2^{\frac{2}{3}}}
 \airy'(t)]n^{-\frac{2}{3}} +
 \end{eqnarray*}
 \begin{eqnarray*}
[ (\frac{5 \alpha -15 c \alpha + 2 c^3 -15 c^2 -56 c -6}{60} +
\frac{c-1}{40} t^3 ) \airy(t)
\end{eqnarray*}
\begin{equation}\label{finallag}
 + \frac{(c-1)(5(c-2)c - 3(2+ 5 \alpha))
 }{60} t \airy'(t)]
n^{-1} + O(n^{-\frac{4 }{3}}) \airy(t) \quad \}
\end{equation}

\begin{figure}[h]\label{fig:CompareL}
\includegraphics{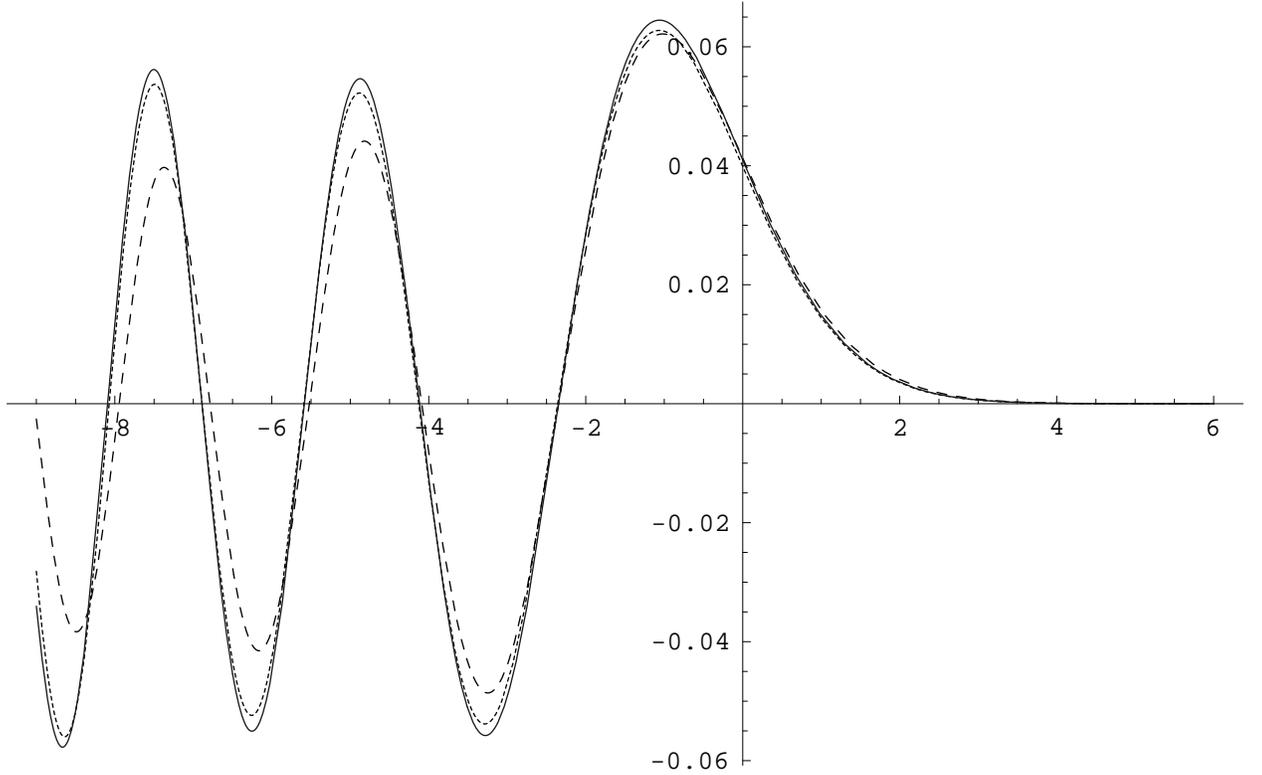}
\caption{For$\; \alpha=-c=1$ and $n=40$, the solid curve
represents $e^{-\frac{\xi^2}{2}} L_{n}^{\alpha}(\xi^{2})$, the
dashed curve is the usual first order approximation of
$e^{-\frac{\xi^2}{2}} L_{n}^{\alpha}(\xi^{2})$ in term of the Airy
function (Which is the first term approximation in
\eqref{finallag}), the doted curve represents our approximation.
These are functions of $t$ where $\xi =(4n +2\alpha
+2c)^{\frac{1}{2}} + \frac{t}{2^{\frac{2}{3}}n^{\frac{1}{6}}} $}
\end{figure}

This is the desired formula for this section. Figure
\ref{fig:CompareL} gives an illustration of our asymptotics.

\clearpage \vspace{3ex} \noindent\textbf{\large Acknowledgements:
} This work was done as part of the author's doctoral thesis
requirements under the supervision of Professor Craig A.~Tracy.
The author would like to thank Professor Tracy for the discussions
that initiated this work and for the invaluable guidance and
support that helped complete it.  Furthermore, the author wishes
to acknowledge helpful discussions with Professor Harold Widom.
This work was supported in parts by the National Science
Foundation under grant DMS--0304414 and by the VIGRE
grant-0135345.


\begin{thebibliography}{10}
\bibitem{Ande1}
G.~Anderson and O. Zeitouni.
\newblock{Lecture Notes On Random Matrices}.
\newblock{Preprint}.

\bibitem{Deif2}
P.~Deift.
\newblock{Orthogonal Polynomials and Random Matrices: A Riemann-Hilbert Approach}.
\newblock {\em American Mathematical Society}. Courant Lecture
Notes 3, 2000.

\bibitem{Deif3}
P.~Deift,
\newblock{Universality for mathematical and physical systems}
\newblock{preprint, arXiv:math-ph/0603038}.

\bibitem{Dien1}
M.~Dieng and C.~A.~Tracy.
\newblock{Application of random matrix theory to multivariate statistics}.
\newblock{preprint, Arxiv:math.PR/0603543}.

\bibitem{Elka1}
N.~El~Karoui.
\newblock {On the largest eigenvalue of Wishart matrices with identity
  covariance when $n$, $p$ and $p/n$ tend to infinity}.
\newblock ArXiv:math.ST/0309355.

\bibitem{Fell2}
W.~Feller.
\newblock {\em An Introduction to Probability Theory and Its
Applications} ,Vol.II.
\newblock Second edition, John Wiley, 1971.

\bibitem{Forr1}
T.~M.~Garoni, P.~J.~Forrester and N.~E.~Frankel.
\newblock{Asymptotic corrections to the eigenvalue density of the
GUE and LUE.}
\newblock arXiv:math-ph/0504053 v1

\bibitem{Gohb1}
I.~Gohberg, S.~Goldberg, and M.~A. Kaashoek.
\newblock {\em {Classes of Linear Operators, Vol. I}}, volume~49 of {\em
  Operator Theory: Advances and Applications}.
\newblock Birkh{\"a}user, 1990.

\bibitem{Gohb2}
I.~Gohberg, S.~Goldberg, and M.~A. Kaashoek.
\newblock {\em {Classes of Linear Operators, Vol. II}}, volume~63 of {\em
  Operator Theory: Advances and Applications}.
\newblock Birkh{\"a}user, 1993.

\bibitem{Gohb3}
I.~C. Gohberg, M.~G. Kre$\breve{i}$n.
\newblock {\em {Introduction to the Theory of Linear Nonselfadjoint Operators}}, volume~18 of {\em
  Translations of Mathematical Monographs}.
\newblock American Mathematical Society, 1969.



\bibitem{Hoch1}
H.~Hochstadt.
\newblock {\em {The Functions of Mathematical Physics}}, volume~23 of {\em
  Pure and Applied Mathematics: A series of texts and Monographs }.
\newblock Wiley-Interscience, 1971.

\bibitem{John1}
I.~M.~Johnstone,
\newblock {On the distribution of the largest eigenvalue in principal component
  analysis},
\newblock {\em Ann. Stats.}, 29(2):295--327, 2001.


\bibitem{Lax1}
P.~D.~Lax.
\newblock{ \em Functional Analysis}
\newblock{ Wiley-Interscience}, 2002.


\bibitem{Meht1}
M.~L.~Mehta.
\newblock {\em Random Matrices, Revised and Enlarged Second Edition}.
\newblock Academic Press, 1991.

\bibitem{Olve1}
F.~W.~J.~Olver.
\newblock{Asymptotics and Special Functions}
\newblock Academic Press, New York, 1974.

\bibitem{Plan1}
M.~Plancherel and W. Rotach.
\newblock{ Sur les valeurs asymptotiques des polynomes d'Hermite}
\newblock Comm. Math. Helv. 1 (1929)227-254.

\bibitem{Sosh1}
A.~Soshnikov.
\newblock {Universality at the Edge of the Spectrum in Wigner Ranom Matrices}.
\newblock {\em J. Stat. Phys.}, 108(5--6):1033--1056, 2002.

\bibitem{Szeg1}
G.~Szeg\"o.
\newblock{Orthogonal Polynomials.}
\newblock American Mathematical Society Colloquium Publications
Volume 23


\bibitem{Trac3}
C.~A.~Tracy and H.~Widom.
\newblock {Level--spacing distributions and the Airy kernel}.
\newblock {\em Commun. Math. Physics}, 159:151--174, 1994.

\bibitem{Trac7}
C.~A.~Tracy and H.~Widom.
\newblock {Fredholm determinants, differential equations and matrix models}.
\newblock {\em Commun. Math. Physics}, 163:33--72, 1994.

\bibitem{Trac2}
C.~A.~Tracy and H.~Widom.
\newblock {On orthogonal and symplectic matrix ensembles}.
\newblock {\em Commun. Math. Physics}, 177:727--754, 1996.

\bibitem{Trac1}
C.~A.~Tracy and H.~Widom.
\newblock {Correlation functions, cluster functions, and spacing distributions
  for random matrices}.
\newblock {\em J. Stat. Phys.}, 92(5--6):809--835, 1998.

\bibitem{Trac4}
C.~A. Tracy and H.~Widom.
\newblock {Airy kernel and Painlev\'e II}.
\newblock In {\em Isomonodromic deformations and applications in physics},
  volume~31 of {\em {CRM Proceedings \& Lecture Notes}}, pages 85--98. Amer.
  Math. Soc., Providence, RI, 2002.

\bibitem{Trac8}
C.~A.~Tracy and H.~Widom.
\newblock {Distribution functions for largest eigenvalues
and their applications}.
\newblock In {\em Proceedings of the International Congress
of Mathematicians, Beijing 2002}, Vol.~I, ed. LI Tatsien,
Higher Education Press, Beijing, pgs.~587--596, 2002.

\bibitem{Trac5}
C.~A.~Tracy and H.~Widom.
\newblock {Matrix kernels for the Gaussian orthogonal and symplectic
  ensembles}.
\newblock {\em Ann. Inst. Fourier, Grenoble}, 55, 2197--2207, 2005.

\bibitem{Whit1}
E.~T.~Whittaker and G.~N.~Watson.
\newblock{\em A Course of Modern Analysis} Fourth Edition
\newblock{Cambridge University Press}, 2004.


\end{thebibliography}
\end{document}